\documentclass{article}

\usepackage[english]{babel}

\usepackage[letterpaper, margin=1in]{geometry}

\usepackage{amsmath}
\usepackage{graphicx}
\usepackage{subfig}
\usepackage[colorlinks=true, allcolors=blue]{hyperref}
\usepackage{authblk}
\usepackage{amssymb}

\usepackage{color, soul}

\providecommand{\keywords}[1]
{
  \small	
  \textbf{\textit{Keywords---}} #1
}

\def\clo#1{\overline{#1}}
\def\text#1{\mbox{#1}}

\def\e#1{ \times 10^{#1}}

\def\AA{\mathcal{A}}
\def\BB{\mathcal{B}}
\def\HH{\mathcal{H}}

\def\KK{\mathcal{K}}

\def\SS{\mathcal{S}}
\def\DD{\mathcal{D}}
\def\Losrc{\Lambda_{\rm osrc}}
\def\EPW{N_{\lambda}}
\def\NG{N_{\rm quad}}

\def\uinc{u_{\rm inc}}
\def\usc{u_{\rm sc}}
\def\utot{u_{\rm tot}}

\graphicspath {{Figures/}}

\title{\textbf{Method of virtual sources using on-surface radiation conditions for the Helmholtz equation}}

\author[1]{Sebastian Acosta}
\author[2]{Tahsin Khajah}

\affil[1]{Department of Pediatrics, Division of Cardiology, Baylor College of Medicine and Texas Children's Hospital}
\affil[2]{Department of Mechanical Engineering, University of Texas at Tyler}

\begin{document}
\maketitle

\begin{abstract}
We develop a novel method of virtual sources to formulate boundary integral equations for exterior wave propagation problems. 
However, by contrast to classical boundary integral formulations, 
we displace the singularity of the Green's function by a small distance $h>0$. As a result, the discretization can be performed on the actual physical boundary with continuous kernels so that any naive quadrature scheme can be used to approximate integral operators. Using on-surface radiation conditions, we combine single- and double-layer potential representations of the solution to arrive at a well-conditioned system upon discretization. The virtual displacement parameter $h$ controls the conditioning of the discrete system. We provide mathematical guidance to choose $h$, in terms of the wavelength and mesh refinements, in order to strike a balance between accuracy and stability. Proof-of-concept implementations are presented, including piecewise linear and isogeometric element formulations in two- and three-dimensional settings. We observe exceptionally well-behaved spectra, and solve the corresponding systems using matrix-free GMRES iterations. The results are compared to analytical solutions for canonical problems. We conclude that the proposed method leads to accurate solutions and good stability for a wide range of wavelengths and mesh refinements. 
\end{abstract}


\section{Introduction}

Methods based on virtual sources are widely employed to approximate the solution to wave radiation and scattering problems governed by the Helmholtz equation. The central idea of these methods is similar to that of integral equations, namely, to represent the solution as a (discrete or continuous) distribution of monopole and/or dipole sources. For classical boundary integral equations, these sources are located on the boundary of the physical domain. For the virtual source methods, the sources are located outside of the physical domain, hence the name ``virtual''. This approach is known by various names including the method of fundamental solutions, generalized multipole techniques, wave superposition method, equivalent source method, and source simulation technique. The similarities and differences between these approaches are explored in \cite{Koopmann1989,Song1991,Mohsen1992,Fairweather2003,Lee2017,Shaaban2023,Martin-Book-2006} and references therein. In recent years, a renewed effort to better characterize the convergence and to expand the applicability of the method of virtual sources has taken place. These efforts are reviewed here \cite{Lee2017,Uscilowska2019,Cheng2020}. Goals include the use of virtual sources as a means to augment the finite element method with Trefftz-type bases \cite{Elsheikh2023, Piltner2019, Zhou2019, Alves2018, Zhou2018}, fast evaluation of integral operators for scattering problems \cite{Bauinger2023, Bauinger2021,Godinho2019, Gaspar2019, Lin2018}, and conditioning/regularization for solving inverse problems \cite{Wang2023,Cheng2020,Ji2019,Karageorghis2018,Geng2015,Karageorghis2011}.

One of the motivations for considering virtual sources is to displace the singularity of the fundamental solution (Green's function) that appears in boundary integral formulations. Once the singularities are placed on a virtual location, then a naive, off-the-shelf quadrature scheme can be employed to approximate the value of the solution in the physical domain including its boundary. However, this idea of displacing these sources can be a double-edged sword. It is precisely the singularity of the Green's function what renders boundary integral equations stably solvable. Hence, the displacement of the sources must be done judiciously. For the method of continuous virtual sources, the displacement distance must be large enough to guarantee the accuracy of the quadrature scheme but simultaneously small enough to control the ill-conditioning of the resulting discrete system of equations.

In this paper, we propose a method of continuous virtual sources that places these sources on a parallel surface shifted by a small distance $h >0$ from the physical boundary of the domain. We represent the solution by a generalized combination of single- and double-layer potentials preconditioned by an on-surface radiation approximation of the Dirichlet-to-Neumann map. We also provide analysis and numerical experiments indicating that $h$ in the order of magnitude of the wavelength $\lambda$ strikes a good balance to regularize the singularity of the fundamental solution as well as to control the conditioning of the resulting matrix. Furthermore, the resulting matrices display exceptionally well-behaved spectra. We also show how the generalized minimal residual (GMRES) method renders accurate solutions in a small number of iterations. The results are compared to analytical solutions for canonical problems to verify accuracy and order of convergence. We present numerical implementations using a simple piecewise linear and a higher order isogeometric (IGA) discretization of the integrals in the layer potentials. Similar to IGA-BEM, no volume parameterization is required \cite{Simpson2014, Chen2020a, Chen2022, Chen2022, Xie2022}. Hence, boundary representations, which are commonly used to generate CAD models, can be adopted directly to alleviate mesh-generation and facilitate shape and topology optimization \cite{Chen2020, Khajah2021, Lian2017, Gao2022}.

\section{Mathematical formulation of the exterior Helmholtz problem}\label{Sec.MathFormulation}

Let us consider a wave propagation problem governed by the Helmholtz equation outside of a bounded domain $\Omega^- \in \mathbb{R}^d$ for dimension $d=2,3$ that represents an impenetrable body with boundary $\Gamma = \partial \Omega^-$. The wave field $u$ under consideration solves the following boundary value problem
\begin{align}
    \Delta u + k^2 u = 0 \qquad & \text{in $\Omega^{+}$,} \label{Eqn.001} \\
    u = f \qquad & \text{on $\Gamma$,} \label{Eqn.002}
\end{align}
subject to the Sommerfeld radiation condition
\begin{align}
    \partial_{r} u - i k u = o\left( r^{\frac{1-d}{2}} \right) , \quad & r=|x| \to \infty \label{Eqn.003}
\end{align}
uniformly for all directions $x/|x|$. Here $\Omega^{+} = \mathbb{R}^d \setminus \clo{\Omega^{-}}$. For simplicity, we focus on the Dirichlet boundary condition \eqref{Eqn.002}, but a similar approach can be followed for a Neumann or a Robin boundary condition. The main underlying assumptions are that the wavenumber $k > 0$ is a constant in $\Omega^{+}$ and that the boundary $\Gamma$ is sufficiently smooth. The well-posedness of the problem \eqref{Eqn.001}-\eqref{Eqn.003} has been established in strong and weak formulations. See for instance \cite{ColtonKressBook2013,NedelecBook2001,McLean2000} and references therein. 

An essential ingredient for our proposed method of virtual sources is Green's representation for the outgoing wave field,
\begin{align}
u(x) = \int_{\Gamma} \left( \partial_{\nu(y)} \Phi(x,y) u(y) - \Phi(x,y) \partial_{\nu} u(y)  \right) d S(y), \qquad x \in \Omega^{+}  \label{Eqn.GreenId}
\end{align}
in terms of its Dirichet data $u$ and Neumann data $\partial_{\nu} u$ on the boundary $\Gamma$. Here $\partial_{\nu}$ represents the derivative in the outward normal direction. The fundamental solution to the Helmholtz equation is
\begin{align}
\Phi(x,y) &= \frac{1}{4 \pi} \frac{e^{i k |x-y|}}{|x-y|}, \quad x \neq y, \quad \text{in $\mathbb{R}^3$, and} \nonumber \\
\Phi(x,y) &= \frac{i}{4} H_{0}^{(1)}(k |x-y|), \quad x \neq y, \quad \text{in $\mathbb{R}^2$.} \nonumber
\end{align}
Since the problem \eqref{Eqn.001}-\eqref{Eqn.003} is uniquely solvable for a given $f$, then the Dirichlet-to-Neumann (DtN) map
\begin{align}
\Lambda: f \mapsto \partial_{\nu}u
\label{Eqn.DtN}
\end{align}
is well-defined in appropriate normed spaces. Consequently, if the DtN map was known a-priori, then \eqref{Eqn.GreenId} would render the solution to 
\eqref{Eqn.001}-\eqref{Eqn.003} in terms of the boundary data $f$ as follows,
\begin{align}
u(x) = \int_{\Gamma} \left( \partial_{\nu(y)} \Phi(x,y) f(y) - \Phi(x,y) \Lambda   f(y) \right) d S(y), \qquad x \in \Omega^{+}. \label{Eqn.GreenId2}
\end{align}
Of course, knowing the action of the DtN map $\Lambda$ amounts to solving \eqref{Eqn.001}-\eqref{Eqn.003} in the first place. Hence, as such, the representations \eqref{Eqn.GreenId} and \eqref{Eqn.GreenId2} offer very little practical advantage. However, useful approximations to the DtN map exist. This is precisely the subject of on-surface radiation conditions (OSRC) developed over the years starting with Kriegsmann et al. \cite{Kriegsmann1987} and followed by Jones \cite{Jones1988,Jones1990,Jones1992}, Ammari \cite{Ammari1998,Ammari1998b}, Calvo et al. \cite{Calvo2004,Calvo2003}, Antoine, Barucq et al. \cite{Antoine1999,Antoine2008,Antoine2001,Antoine2006,
Barucq2003,Barucq2010a,Barucq2012,Alzubaidi2016} and Darbas, Chaillat, Le Lou\"er et al. \cite{Antoine2004,Antoine2005,Antoine2006,Antoine2007,Darbas2006,Darbas2013,
Darbas2015,Chaillat2015,Chaillat2017}. See also \cite{Atle2007,Acosta2015f,Acosta2017c,Acosta2021d,Stupfel1994,Murch1993,Teymur1996,Yilmaz2007,Medvinsky2010,Chniti2016,Alzubaidi2016} among others. 

One of the applications of OSRC is to provide a well-conditioned formulation for the representation of the solution $u$ as generalized combination of single- and double- layer potentials and the corresponding boundary integral equations. See the work of Antoine, Darbas, Chaillat, Le Lou\"er \cite{Antoine2004,Antoine2005,Antoine2007,Darbas2006,Darbas2013,Darbas2015,Chaillat2015,Chaillat2017,Chaillat2021} including a recent review \cite{Antoine2021}. The \textit{single-layer} potential $\SS$ and \textit{double-layer} potential $\DD$ are defined by 
\begin{align}
(\mathcal{S} v )(x) = \int_{\Gamma} v(y) \Phi(x,y) dS(y) \qquad \text{and} \qquad
(\mathcal{D} v )(x) = \int_{\Gamma} v(y) \partial_{\nu(y)} \Phi(x,y) dS(y), \qquad x \in \Omega^{+}. \label{Eqn.SingleDoubleLP}
\end{align}
Hence, if a practical approximation $\Losrc$ of the DtN map $\Lambda$ is available, then \eqref{Eqn.GreenId2} suggests that the solution $u$ should be represented as 
\begin{align}
u(x) = \left( (\mathcal{D} - \mathcal{S} \Losrc) v \right) (x)
\qquad x \in \Omega^{+}, \label{Eqn.CombinedLP}
\end{align}
for some unknown density $v$. The reason why \eqref{Eqn.CombinedLP} leads to a well-conditioned boundary integral equation is that the boundary trace of $(\mathcal{D} - \mathcal{S} \Losrc)$ approximates the identity operator as suggested by \eqref{Eqn.GreenId2}, as long as $\Losrc$ is a good approximation to $\Lambda$. In the next section, we use these insights to formulate a well-conditioned method of virtual sources.

\section{Method of virtual sources using OSRC} \label{Sec.MVCS}

Here we pose the method of virtual sources at the continuous level. The starting point is to choose a shifted surface contained in $\Omega^{-}$ where the virtual sources are placed. Specifically, we consider a parallel surface described by
\begin{align}
\Gamma_{h} = \left\{ z = y - h \nu(y) ~:~ y \in \Gamma \right\} \label{Eqn.ParallelSurface}
\end{align} 
with a parameter $h>0$, and where $\nu(y)$ is the outward unit normal vector at a point $y \in \Gamma$. When $\Gamma$ is sufficiently smooth and $h$ sufficiently small, then $\Gamma_{h}$ is well defined. Moreover, the normal vector $\nu_{h}(z)$ of the parallel surface $\Gamma_{h}$ coincides with the normal vector $\nu(y)$ of $\Gamma$ for all $y \in \Gamma$. See details in \cite[\S 6.2]{Kress-Book-1999}. 

Now we seek the solution to the boundary value problem \eqref{Eqn.001}-\eqref{Eqn.003} in the form of a generalized combined layer potential emanating from the ``virtual" surface $\Gamma_{h}$,
\begin{align}
u(x) &= \int_{\Gamma_{h}} \left( \partial_{\nu_{h}(z)} \Phi(x,z) - \Phi(x,z) \Losrc  \right) v_{h}(z)  d S_{h}(z) \nonumber \\
&= \int_{\Gamma}  \left( \partial_{\nu(y)} \Phi(x,y-h\nu(y)) - \Phi(x,y-h\nu(y)) \Losrc  \right) v(y) \left( 1 - 2 h \HH(y) + h^2 \KK(y) \right)  d S(y).
\label{Eqn.VirtualSource}
\end{align}
The second equality above is simply due to a change of variables $z = y - h \nu(y)$, defining $v(y) = v_{h}(z)$, and the fact that surface elements on $\Gamma$ and $\Gamma_{h}$ are related by
\begin{align}
dS_{h}(z) = \left( 1 - 2 h \HH(y) + h^2 \KK(y) \right)  d S(y)
\label{Eqn.SurfaceElement}
\end{align}
where $\HH$ and $\KK$ are the mean and Gauss curvatures of $\Gamma$, respectively.

By virtue of the fundamental solution $\Phi$, the potential in \eqref{Eqn.VirtualSource} solves the Helmholtz equation in $\Omega^{+}$ and the Sommerfeld radiation condition at infinity. The potential also satisfies the Dirichlet boundary condition \eqref{Eqn.002} provided that the density $v$ is a solution to the following integral equation 
\begin{align}
\int_{\Gamma}  \left( \partial_{\nu(y)} \Phi(x,y-h\nu(y)) - \Phi(x,y-h\nu(y)) \Losrc  \right) v(y) \left( 1 - 2 h \HH(y) + h^2 \KK(y) \right)  d S(y) = f(x), \quad x \in \Gamma.
\label{Eqn.MainIntEqn}
\end{align}
Note that the factor $\left( 1 - 2 h \HH + h^2 \KK \right)$ could be absorbed into the unknown density $v$ so that there is no need to compute curvatures at this point.

Introduce an integral operator $\AA : H^{s}(\Gamma) \to H^{s}(\Gamma)$, $s\geq 1/2$, by
\begin{align}
(\AA v)(x) = \int_{\Gamma}  \left( \partial_{\nu(y)} \Phi(x,y-h\nu(y)) - \Phi(x,y-h\nu(y)) \Losrc  \right) v(y) \left( 1 - 2 h \HH(y) + h^2 \KK(y) \right)  d S(y).
\label{Eqn.IntOp}
\end{align}
Then we can express the integral equation \eqref{Eqn.MainIntEqn} in the form
\begin{align}
\AA v = f. \label{Eqn.100}
\end{align}

One of the challenges of the method of virtual sources is the ill-conditioning of the governing operator and its discrete versions. Since the integral operator $\AA$ has a smooth kernel for any $h>0$, then $\AA$ is a compact operator in any Sobolev degree and it cannot have a bounded inverse. In fact, zero belongs to the spectrum of $\AA$, that is, $\AA$ is not injective or otherwise the eigenvalues of $\AA$ accumulate at zero. This is the essential drawback of displacing the singularity of the fundamental solution as proposed by the method of virtual sources.

We proceed to construct an approximate inverse for the operator $\AA$ based on the properties of the Green's identity discussed in Section \ref{Sec.MathFormulation}. Note that if we replace $\Losrc$ with the exact DtN map $\Lambda$ in \eqref{Eqn.VirtualSource} and in the definition \eqref{Eqn.IntOp} of the operator $\AA$, then this version of the operator $\AA$ would simply propagate the density $v$ from the virtual surface $\Gamma_{h}$ to the physical surface $\Gamma$. Such outgoing wave field solves $\partial_{\nu} v = \Lambda v$ which means that
\begin{align}
\AA v \approx e^{h \Lambda} v. \label{Eqn.102}
\end{align}
Hence, an approximate inverse to the operator $\AA$ can be defined as follows
\begin{align}
\BB = I - h \Losrc \approx e^{-h \Losrc} \approx e^{-h \Lambda}. \label{Eqn.104}
\end{align}
We can expect $\BB$ to be a good preconditioner for $\AA$ as long as $h$ is sufficiently small and the on-surface radiation operator $\Losrc$ is a good approximation to the exact DtN map $\Lambda$. Hence, instead of solving \eqref{Eqn.100}, we propose to solve the following equation as our well-conditioned method of virtual sources,
\begin{align}
\BB \AA v = \BB f. \label{Eqn.106}
\end{align}

\section{Piecewise linear discretization in 2D} \label{Section.LinearDiscrete}

As a first proof-of-concept for the proposed method of virtual sources with OSRC, we implemented a discretization using piecewise linear elements. The proposed method takes the form of the integral equation \eqref{Eqn.106} with continuous kernel provided $h>0$. So here we apply a Nystrom method for continuous kernels by a sequence of convergent quadrature operators  
\begin{align}
(\AA_{N} v)(x) &= \sum_{n=1}^{N} \alpha_{n,N} \bigg[ K_{D}(x,y_{n,N} - h\nu(y_{n,N})) v(y_{n,N}) \nonumber \\ & \quad - K_{S}(x,y_{n,N} - h\nu(y_{n,N})) v(y_{n,N}) (\Losrc v)(y_{n,N}) \bigg] , \qquad x \in \Gamma \cup \Omega^{+},
\end{align}
with a corresponding sequence of quadrature points $\{ y_{n,N} \} \in \Gamma$, and quadrature weights $\{ \alpha_{n,N} \}$. Here we assume that the quadrature points are distributed near uniformly along the boundary $\Gamma$ so that the size of the elements is also nearly uniform, and the number of elements per wavelength $\EPW \sim N \lambda / |\Gamma|$.
Here the continuous kernels $K_{D}(x,y) = \partial_{\nu(y)} \Phi(x,y)$ and $K_{S}(x,y) = \Phi(x,y)$. The approximate DtN map $\Losrc$ acts on the piecewise linear function $v$ through the Pad\'e form described in the Appendix \ref{Section.Pade} using a piecewise linear finite element scheme for the tangential derivatives along $\Gamma$, and an angle of rotation $\theta = \pi/2$.

The solution $v$ of \eqref{Eqn.106} is approximated by the piecewise linear solution $v_{N}$ of the following $N \times N$ linear system
\begin{align}
(\BB_{N} \AA_{N} v_{N})(x_{n,N}) = (\BB_{N} f)(x_{n,N}), \qquad n=1, 2, ..., N, \label{Eqn.LinearSystem01}
\end{align}
where the collocation points $\{ x_{n,N} \}$ coincide with the quadrature points $\{ y_{n,N} \}$. We propose to solve this system by a matrix-free GMRES method. 

\begin{figure}[ht]
\centering
\captionsetup{font=small}
\subfloat[Numerical solution]{
\includegraphics[height=0.28 \textwidth, trim = 0 0 0 0, clip]{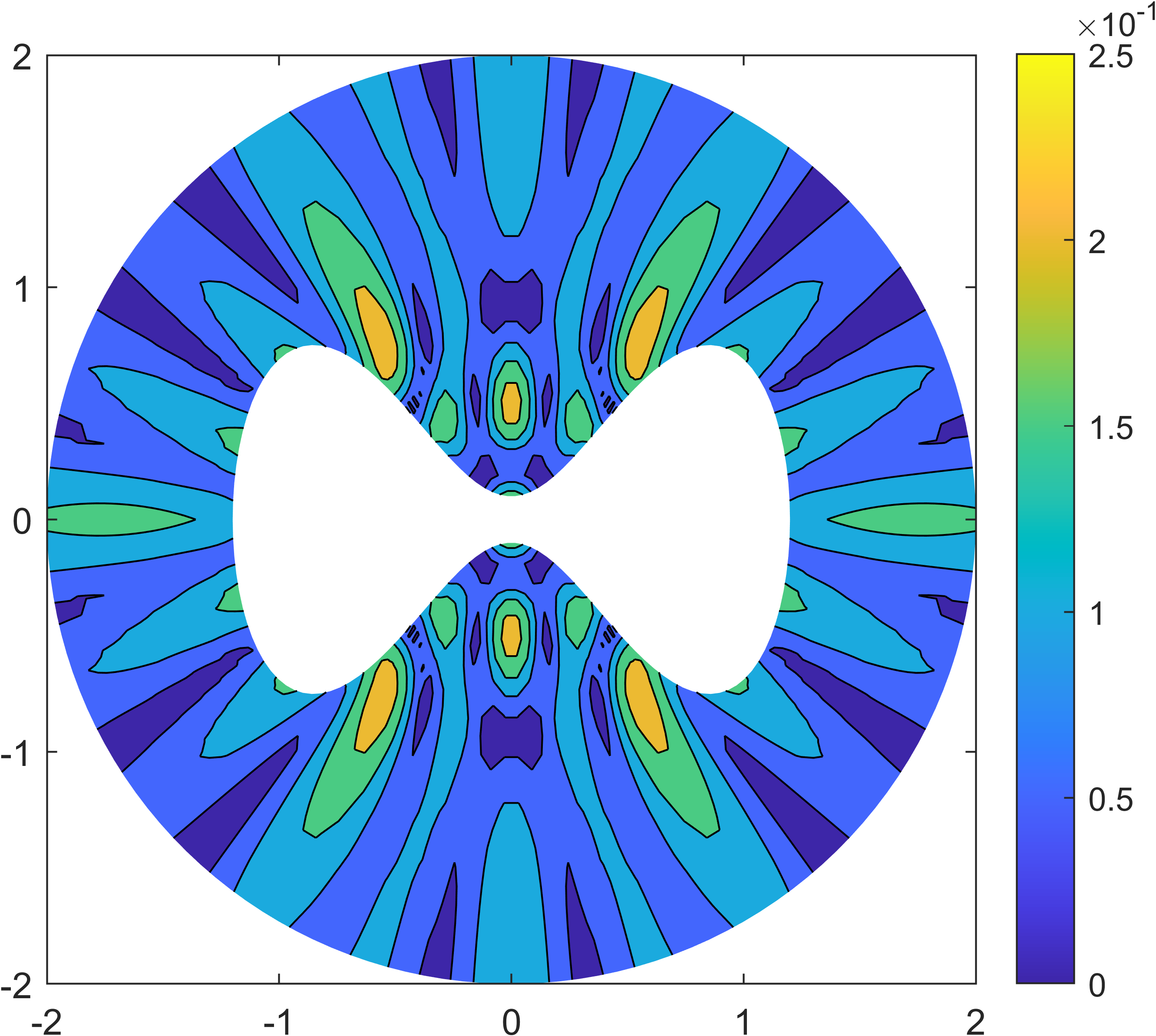}
\label{subfig:NumSol}
}
\subfloat[Error ($\log_{10}$ scale)]{
\includegraphics[height=0.28 \textwidth, trim = 0 0 0 0, clip]{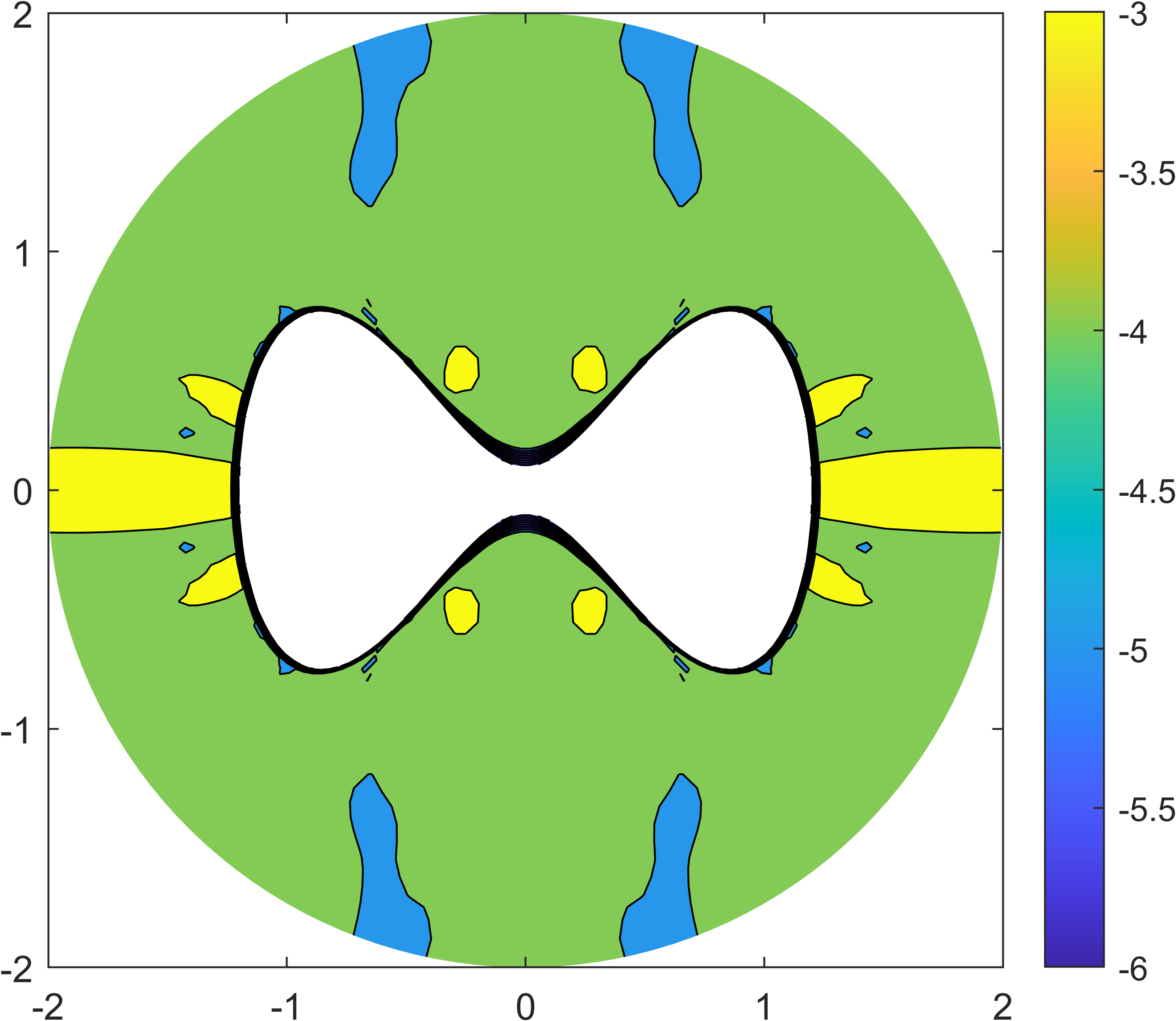}
\label{subfig:Error}
}
\subfloat[GMRES residual]{
\includegraphics[height=0.28 \textwidth, trim = 0 0 0 0, clip]{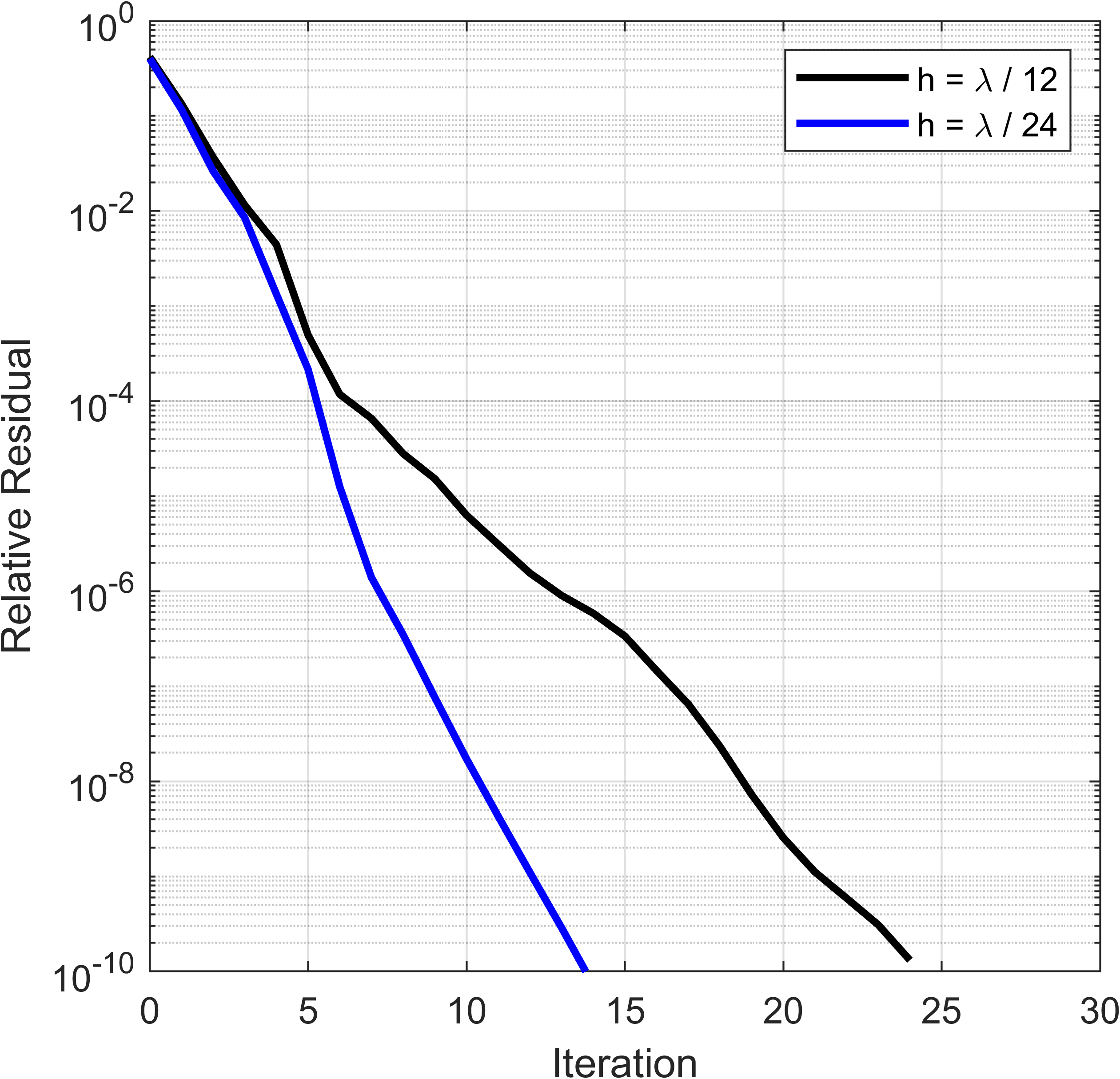}
\label{subfig:GMRES_res}
}
\caption{(a) Numerical solution for $k=4\pi$, (b) the error ($\log_{10}$ scale) in a neighborhood of the boundary $\Gamma$, and (c) the residual for the GMRES iterations for the cases $h=\lambda/12$ and $h=\lambda/24$. The boundary $\Gamma$ is described by \eqref{Eqn.GammaParam}. The numerical solution was obtained using the proposed method of virtual sources with $\EPW=12$, a Pad\'e approximation $\Losrc$ of the DtN map using $4$ terms.}
\label{Fig.ExactSol_and_comparison}
\end{figure}

For a two-dimensional numerical experiment, we let $\Gamma$ be a smooth curve parameterized as follows,
\begin{align}
x(t) = (1 + 0.2 \cos(2 t)) \cos(t) \quad \text{and} \quad
y(t) = (1 + 0.9 \cos(2 t)) \sin(t) \label{Eqn.GammaParam}
\end{align}
for $t \in [0, 2 \pi)$. This shape is illustrated in Figure \ref{Fig.ExactSol_and_comparison}. We manufactured a known exact solution to the boundary value problem \eqref{Eqn.001}-\eqref{Eqn.003} with this non-canonical boundary $\Gamma$ by setting $f$ equal the boundary trace of the following radiating field,
\begin{align}
F(x) = \sum_{j=1}^{4} \Phi(x,y_{j}) \label{Eqn.Field}
\end{align}
where the four point sources $\{ y_{j} \}$ are located inside $\Gamma$ so that $F$ is a radiating solution to the Helmholtz equation in $\Omega^{+}$. Specifically, $y_{1} = (3/4, -1/2)$, $y_{2} = (3/4, 1/2)$, $y_{3} = (-3/4, 1/2)$, and $y_{4} = (-3/4, -1/2)$. As a consequence, the exact solution to \eqref{Eqn.001}-\eqref{Eqn.003} is $u = F$ in $\Omega^{+}$. Having the exact solution available allows us to compute a discrete $l^2$ relative error,
\begin{align}
\text{Relative Error} = \frac{ \sum_{j=1}^{J} | 
(\AA_{N} v)(z_{j}) - u(z_{j}) |^2 }{ \sum_{j=1}^{J} |u(z_{j})|^2 }  \label{Eqn.RelError}
\end{align}
for a collection of points $\{ z_{j} \}_{j=1}^{J} \subset \Omega^{+}$, in the vicinity of the boundary $\Gamma$. See Figure \ref{Fig.ExactSol_and_comparison} for an illustration of the numerical solution and the pointwise error in the vicinity of $\Gamma$, and the progression of the residuals from the GMRES iterations for the cases $h = \lambda/12$ and $h=\lambda/24$. We note that the GMRES iterations converge faster for $h=\lambda/24$ than for $h=\lambda/12$ as expected from the spectral characterzation of this method analyzed in the Appendix \ref{Section.Spectral}.

\begin{figure}[ht]
\centering
\captionsetup{font=small}
\subfloat[$k=4 \pi$]{
\includegraphics[height=0.24 \textwidth, trim = 2 0 2 0, clip]{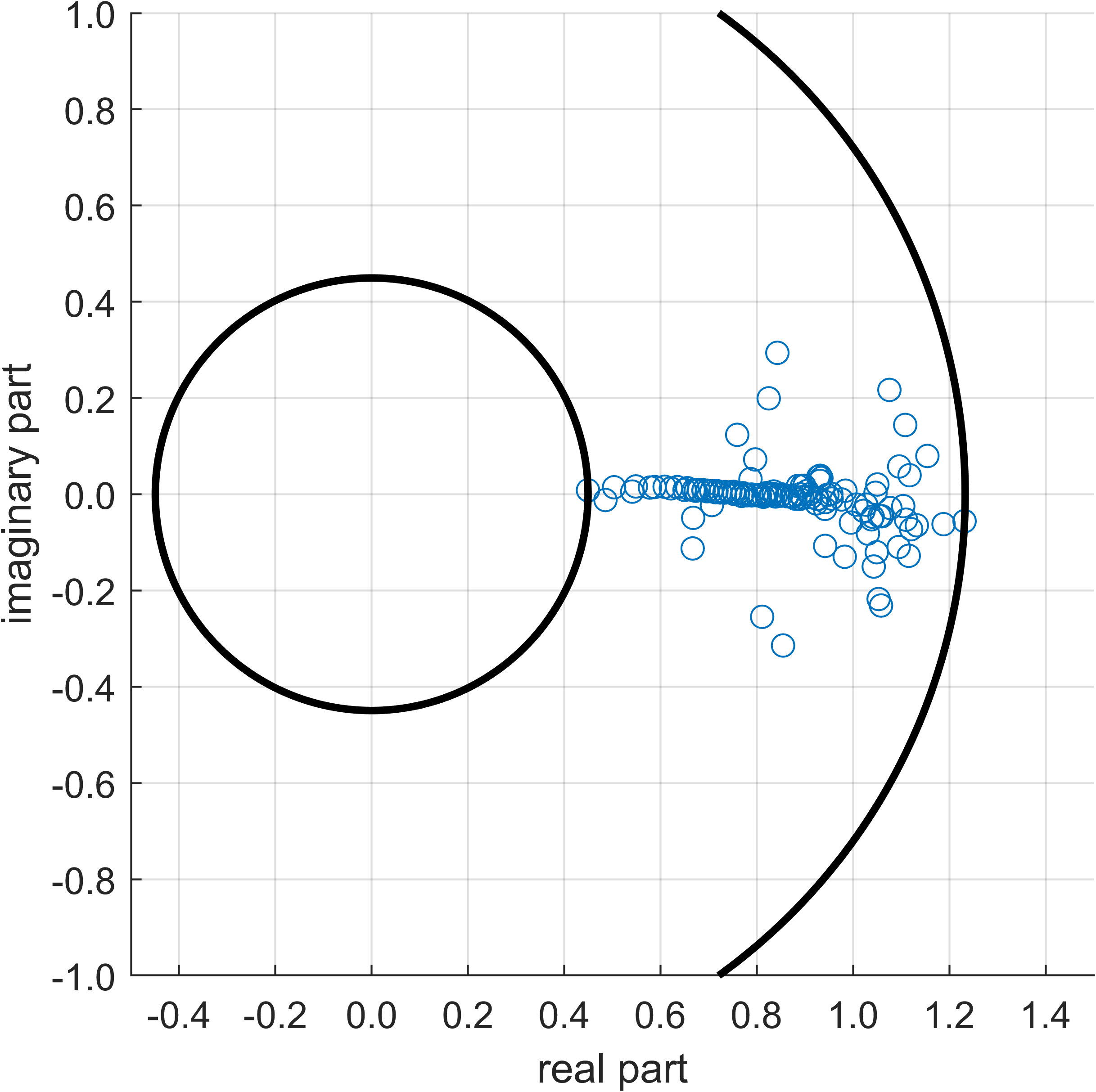}
\label{subfig:K4pi}
}
\subfloat[$k=8 \pi$]{
\includegraphics[height=0.24 \textwidth, trim = 2 0 2 0, clip]{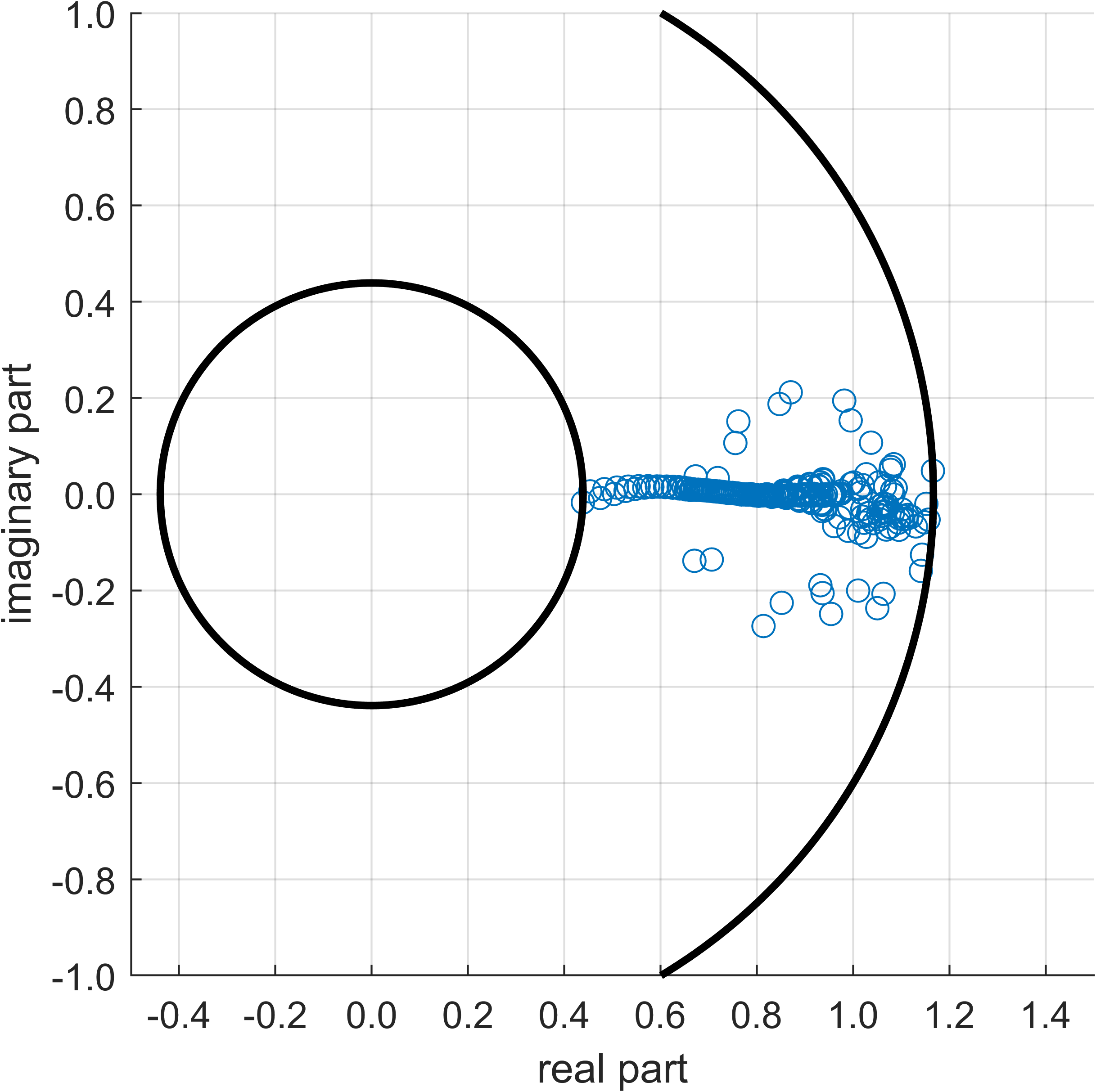}
\label{subfig:K8pi}
}
\subfloat[$k=16 \pi$]{
\includegraphics[height=0.24 \textwidth, trim = 2 0 2 0, clip]{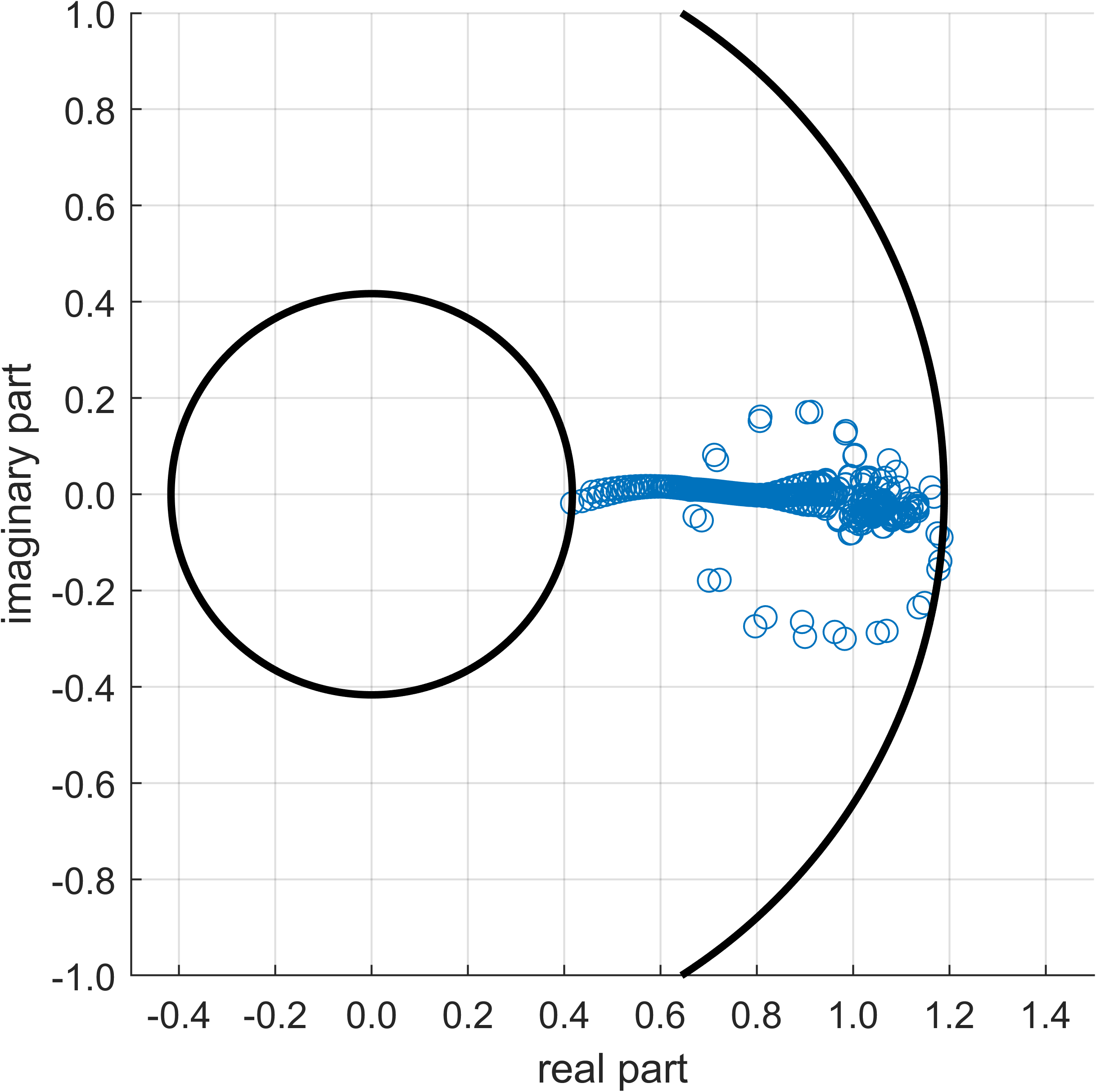}
\label{subfig:K16pi}
} 
\subfloat[$k=32 \pi$]{
\includegraphics[height=0.24 \textwidth, trim = 2 0 2 0, clip]{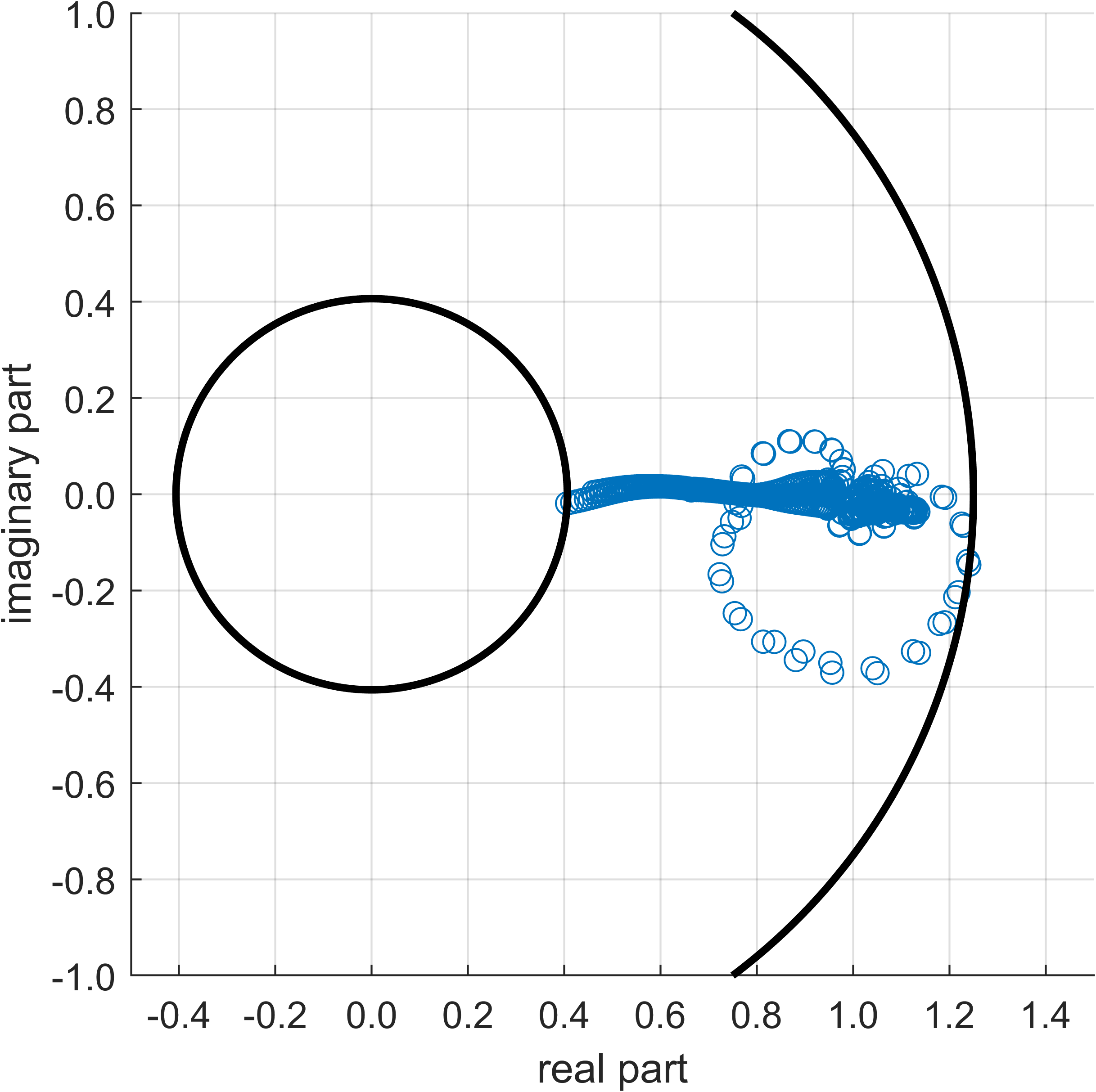}
\label{subfig:K32pi}
}
\caption{Spectrum of the preconditioned matrix $\BB_{N}\AA_{N}$ for wavenumbers $k=4 \pi$, $k=8 \pi$, $k=16 \pi$, and $k=32 \pi$. The circles's radii correspond to the largest and smallest magnitude of the spectrum. In all cases, the discretization was done using $\EPW = 12$ along the boundary $\Gamma$. The corresponding degrees of freedom are $N=150$, $N=302$, $N=604$, and $N=1206$, respectively. The virtual source displacement distance is given by $h = \lambda / \EPW$ where $\lambda = 2\pi / k$ is the wavelength.}
\label{fig:SpectrumIncreasingFreq}
\end{figure}

For this numerical setup, the spectrum of the proposed pre-conditioned matrix $(\BB_{N} \AA_{N})$ is shown in Figure \ref{fig:SpectrumIncreasingFreq} for increasing wavenumber $k=4\pi$, $8\pi$, $16\pi$, and $32\pi$ but fixed number of elements per wavelength $\EPW=12$. The displacement distance for the virtual sources is $h=\lambda/12$ where $\lambda=2\pi/k$ is the wavelength of the fields. 
We observe that the upper- and lower-bounds of the spectrum are practically independent of the frequency. As shown in the Appendix \ref{Section.GaussQuadPhi}, we expect the numerical quadrature to remain accurate near the shifted singularity, and simultaneously the matrix $(\BB_{N} \AA_{N})$ to remain well-conditioned for this choice of $h$ as the wavenumber $k$ increases, as long as $\EPW$ is kept constant. This is a tremendous advantage over the conventional method of virtual sources where the eigenvalues approach zero exponentially fast as $k$ increases and $\EPW$ is kept constant because the number of degrees of freedom $N$ would also increase proportionally to $k$. See \cite{Kress1986} and the discussion in the Appendix \ref{Section.Spectral}.

We numerically explore the spectral behavior and accuracy for constant wavenumber $k=4\pi$, and increasing number of elements per wavelength $\EPW=12$, $\EPW=24$, $\EPW = 48$ and $\EPW = 96$, and for $h= \lambda/ \EPW^{\beta}$. As discussed in the Appendix \ref{Section.Spectral}, a balance between accuracy and stability can be achieved by setting $h= \lambda/ \EPW^{\beta}$ for some $0 \leq \beta \leq 1$. Table \ref{tab:table1} displays the spectral condition number of $(\BB_N\AA_N)$, defined as 
\begin{align}
\text{cond}(A)  = \frac{ \max_{\lambda_{n} \in \sigma(A)} |\lambda_n| }{\min_{\lambda_{n} \in \sigma(A)} |\lambda_n|}\label{Eqn.CondNumber}
\end{align}
where $\sigma(A)$ denotes the set of eigenvalues of a generic matrix $A$, for mesh refinements and for the cases $\beta=0$, $\beta=1/2$, and $\beta=1$. Note that $\beta=0$ induces an exponential increase in condition number as the mesh is refined. The governing matrix is still able to be inverted accurately leading to small relative error in the numerical solution. However, the drastic increase in condition number is worrisome and it can be catastrophic for noisy data. For $\beta=1$, the condition number is controlled but the scheme suffers from lack of accuracy due to the small displacement of the Green function's singularity which limits the accuracy of numerical quadrature. See Appendix \ref{Section.GaussQuadPhi} for details.  
For $\beta=1/2$, a compromise between stability and accuracy is reached, where the condition number still increases but not as drastically as in the case $\beta=0$, and the relative error decreases as the mesh is refined. Notice that the order of convergence with mesh refinement is at least quadratic.

\begin{table}
\centering
\captionsetup{font=small}
\caption {\label{tab:table1} Condition number \eqref{Eqn.CondNumber} and relative error \eqref{Eqn.RelError} for refinements of the mesh $\EPW = 12$, $24$, $48$ and $96$, and various choices for the parameter $\beta$ to govern the virtual sources displacement distance $h \sim \lambda / \EPW^{\beta}$. The wavenumber $k=4 \pi$ is constant and $\lambda = 2 \pi / k$ is the wavelength. The number of Gauss nodes per element $\NG=4$ is also fixed.} 
\begin{tabular}{lllllllll}
\hline
 & \multicolumn{2}{c}{$\beta=0$} & & \multicolumn{2}{c}{$\beta=1/2$} & &  \multicolumn{2}{c}{$\beta=1$}  \\ 
\cline{2-3} \cline{5-6} \cline{8-9} 
$\EPW$ &  Cond. Number  & Rel. Error  &  & Cond. Number  & Rel. Error &  &  Cond. Number  & Rel. Error  \\  

$12$   &  $1.11 \times 10^{1}$  & $5.02 \times 10^{-3}$  & & $1.11 \times 10^{1}$  & $5.02 \times 10^{-3}$ & & $1.11 \times 10^{1}$  & $5.02 \times 10^{-3}$  \\

$24$   &  $2.41 \times 10^{2}$  & $1.96 \times 10^{-4}$  & & $4.47 \times 10^{1}$ & $1.10 \times 10^{-3}$   & & $1.34 \times 10^{1}$  & $4.50 \times 10^{-3}$  \\
  
$48$   &  $3.89 \times 10^{5}$  & $1.35 \times 10^{-7}$ & & $7.02 \times 10^{2}$ & $2.05 \times 10^{-4}$  & & $2.96 \times 10^{1}$  & $4.39 \times 10^{-3}$  \\

$96$   &  $4.37 \times 10^{11}$  & $3.58 \times 10^{-12}$ & & $2.60 \times 10^{4}$ & $1.48 \times 10^{-5}$  & & $4.97 \times 10^{1}$  & $4.39 \times 10^{-3}$  \\
  \hline
\end{tabular}
\end{table}

\section{Isogeometric discretization} \label{Section.IGADiscrete}

Here we describe the implementation of the proposed method of virtual sources using an isogeometric analysis (IGA) formulation. As it is well-known \cite{Hughes2005}, the IGA formulation employs the same basis function representation for the geometry of the domain and for the finite element space to approximate the solution of the wave problem. Specifically, non-uniform rational B-splines (NURBS) are used to approximate the domain and the solution to the integral equation \eqref{Eqn.106}, including the implementation of the OSRC approximation $\Losrc$ of the DtN map. It has been shown that IGA formulations can yield higher accuracy per degree of freedom compared to conventional finite element methods, and reduced pollution error for wave propagation problems \cite{Hughes2005,Khajah2019b,Dsouza2021,Shaaban2023}. Combination of IGA and high-order local absorbing boundary conditions led to very accurate scattering analyses \cite{Khajah2019,Dsouza2021,Atroshchenko2022}. Our implementation of the IGA method follows closely \cite{Antoine2022,Khajah2023} for the definition of the B-splines, the NURBS, and for the calculation of surface curvatures.

\subsection{Manufactured radiating solution for a non-canonical 2D boundary} \label{Sec.SubsectionManuf_IGA}

The setup for the numerical experiment is the same as described in Section \ref{Section.LinearDiscrete} with the boundary $\Gamma$ given by \eqref{Eqn.GammaParam} and Dirichlet boundary data $f$ given by the boundary trace of \eqref{Eqn.Field} which is the exact solution to the boundary value problem \eqref{Eqn.001}-\eqref{Eqn.003} is $u = F$ in $\Omega^{+}$. For the IGA implementation, let $\EPW$ denote the number of elements per wavelength, and $\NG$ the number of Gauss nodes per element to approximate the integrals appearing in the single- and double-layer potentials with displaced sources. In all cases, we use B-splines of degree 3 as a basis to represent the boundary $\Gamma$ and the solution to the wave problem, and Pad\'e approximation $\Losrc$ of the DtN map using $4$ terms. The displacement distance for the virtual sources is again denoted by $h>0$. The IGA discretized version of the governing equation \eqref{Eqn.106} was solved using the GMRES iterative method.

\begin{figure}[h]
\centering
\captionsetup{font=small}
\subfloat[Numerical solution]{
\includegraphics[height=0.28 \textwidth, trim = 0 0 0 0, clip]{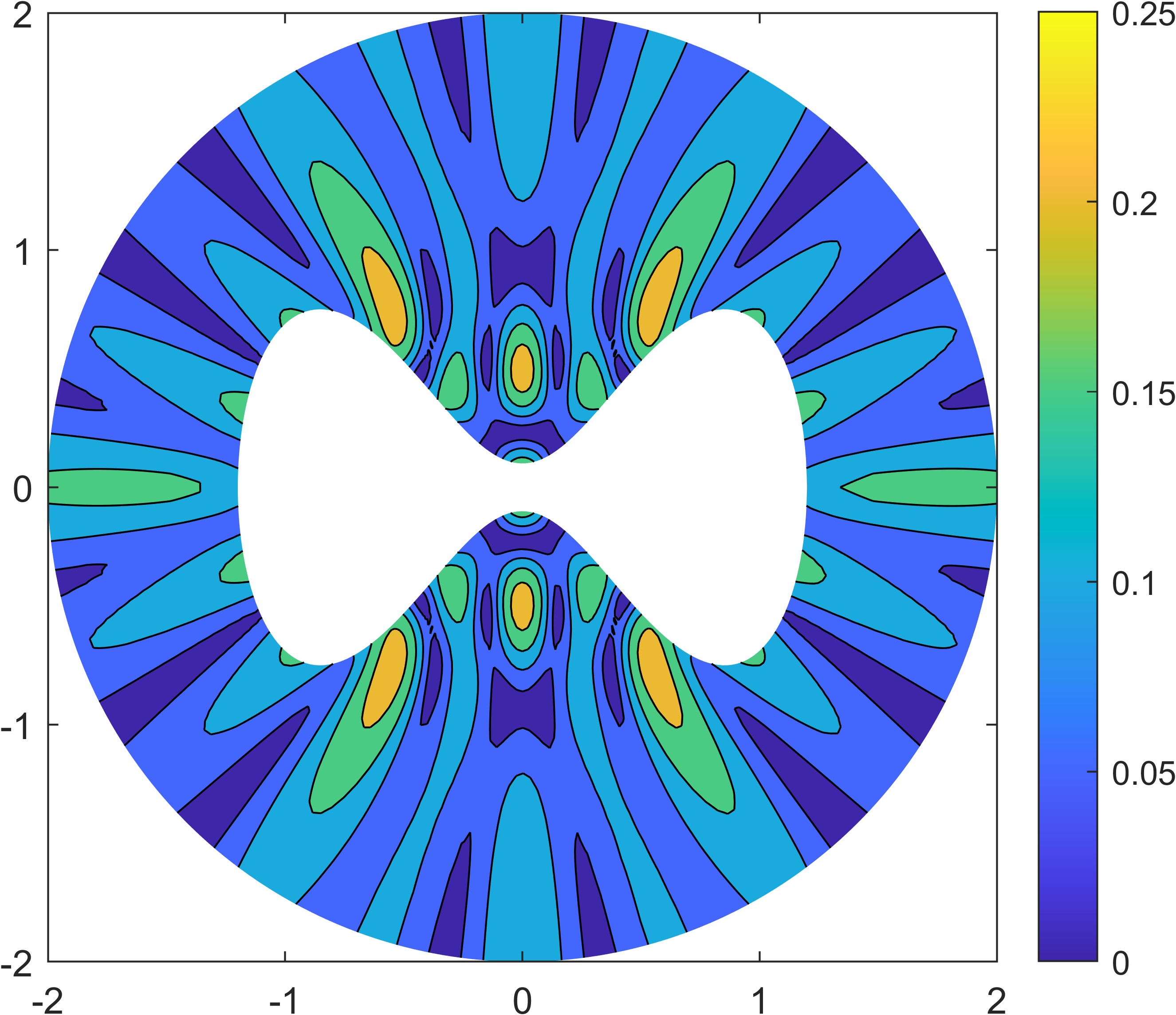}
\label{subfig:NumSol2}
}
\subfloat[Error ($\log_{10}$ scale) ]{
\includegraphics[height=0.28 \textwidth, trim = 0 0 0 0, clip]{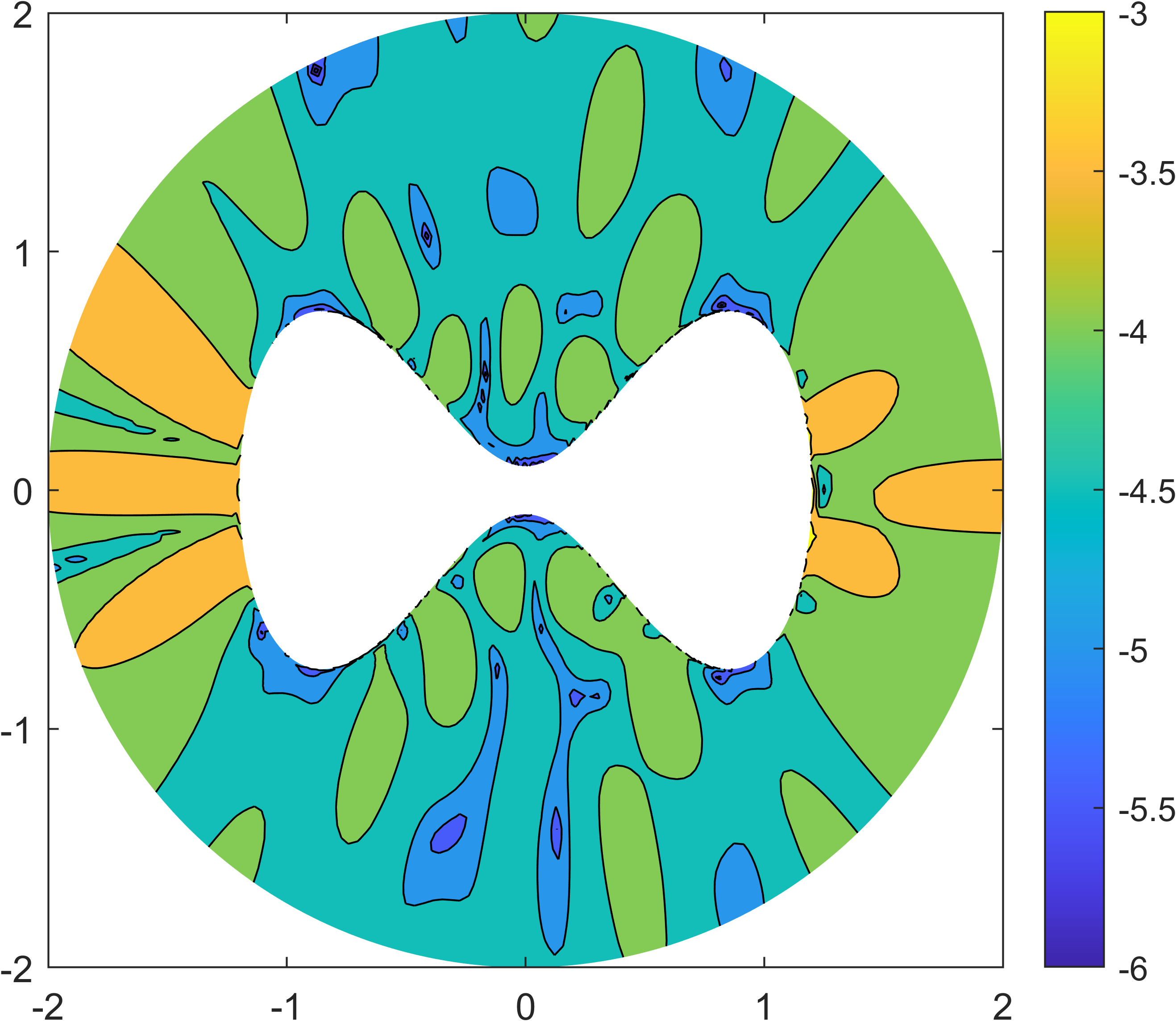}
\label{subfig:Error2}
}
\subfloat[GMRES residual]{
\includegraphics[height=0.28 \textwidth, trim = 0 0 0 0, clip]{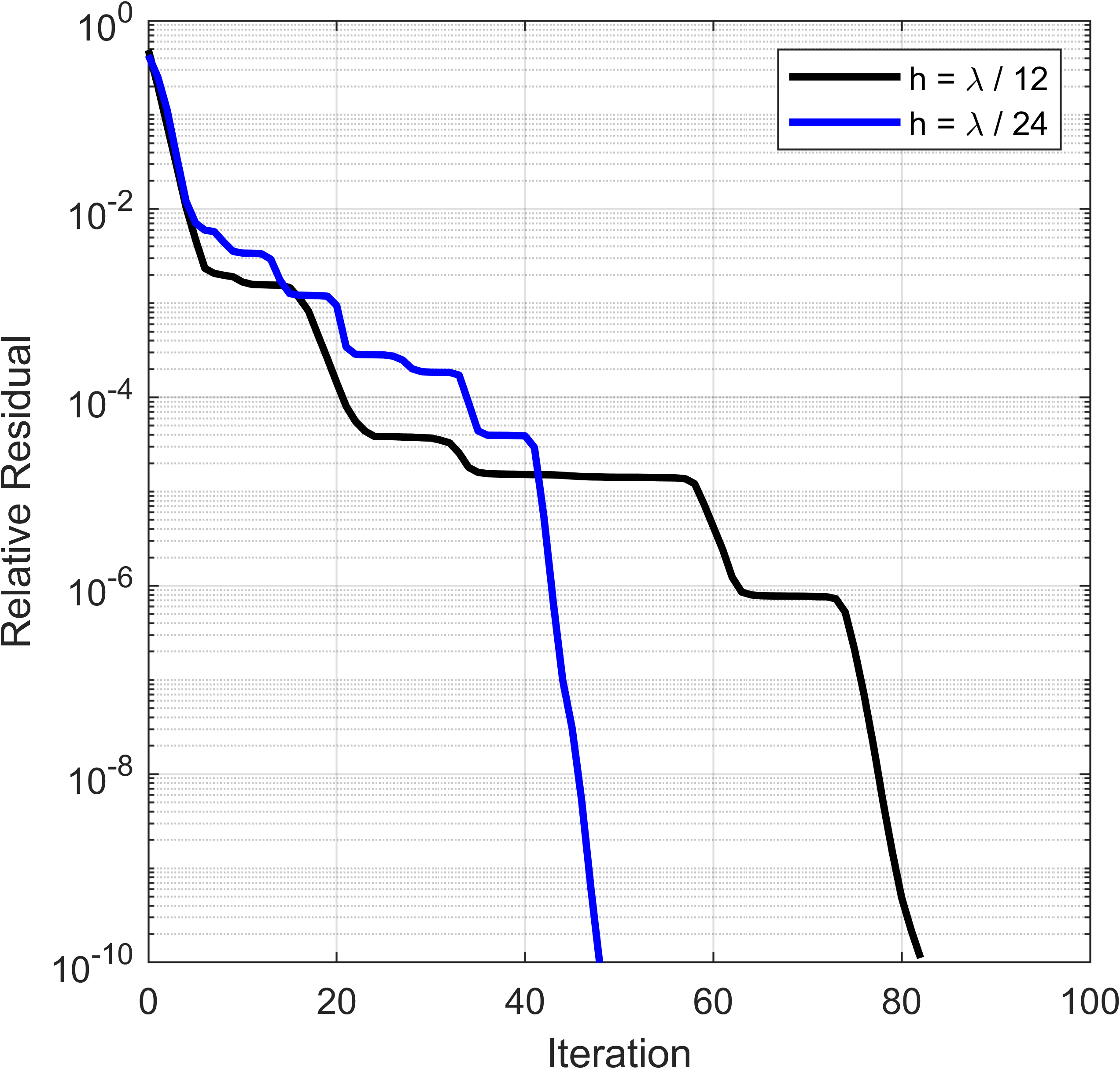}
\label{subfig:GMRES_res2}
}
\caption{ (a) Numerical solution using the IGA implementation for $k=4\pi$, (b) the error ($\log_{10}$ scale) in a neighborhood of the boundary $\Gamma$, and (c) the residual for the GMRES iterations for the cases $h = \lambda/12$ and $h = \lambda/24$. The boundary $\Gamma$ is described by \eqref{Eqn.GammaParam}. The numerical solution was obtained using the proposed method of virtual sources with $\EPW=12$, a Pad\'e approximation $\Losrc$ of the DtN map using $4$ terms, and Gaussian quadrature with $\NG = 4$ nodes per elements.}
\label{Fig.IGA_ExactSol_and_comparison}
\end{figure}

Our first objective is to fix the virtual source displacement distance at $h = \lambda/24$ and explore the error in terms of the number of elements per wave length $\EPW$ and number of Gauss nodes per element $\NG$. This process is illustrated in Table \ref{tab:Error1} for wavenumber $k=4 \pi$ fixed. We observe that by increasing the number of Gauss nodes, the error can be made smaller by several orders of magnitude. At some point, however, the error stagnates. This is expected as the error from numerical integration has been made sufficiently small for the error from the finite-dimensional NURBS basis to dominate. Hence, by then increasing the number of elements per wavelength $\EPW$ we decrease the total error further. As an illustrative example, the numerical solution, error profile and evolution of the residual over the GMRES iterations are illustrated in Figure \ref{Fig.IGA_ExactSol_and_comparison}. Compared to the solution using piecewise linear elements from Section \ref{Section.LinearDiscrete} illustrated in Figure \ref{Fig.ExactSol_and_comparison}, the error is more than one order of magnitude smaller.

\begin{table}[h]
\centering
\captionsetup{font=small}
\caption { \label{tab:Error1} Relative error for the radiating problem described in subsection \ref{Sec.SubsectionManuf_IGA} for various values of the number of elements per wavelength $\EPW$ and number of Gauss nodes per element $\NG$.} 
\begin{tabular}{cccc}
\hline
$\NG$ & $\EPW=6$ & $\EPW=12$ &  $\EPW=24$ \\  

$4$ & $ 3.12\e{-2} $ & $ 1.54\e{-3} $ & $ 1.64\e{-5} $  \\

$5$ & $ 8.40\e{-2} $ & $ 4.80\e{-4} $ & $ 1.42\e{-6} $  \\
  
$6$ & $ 1.00\e{-2} $ & $ 1.47\e{-4} $ & $ 1.88\e{-7} $  \\

$7$ & $ 3.35\e{-3} $ & $ 4.45\e{-5} $ & $ 2.54\e{-8} $  \\

$8$ & $ 2.18\e{-3} $ & $ 1.36\e{-5} $ & $ 5.69\e{-9} $  \\

\hline
\end{tabular}
\end{table}

The second objective is to explore the conditioning of the linear system in terms of the mesh refinement and the virtual source displacement distance. Such a refinement is illustrated in Table \ref{tab:table2} for various choices for the virtual source displacement progression $h \sim \lambda / \EPW^{\beta}$ for $\beta=0$, $\beta=1.5$ and $\beta=1$. As in the case of piecewise linear elements, for $\beta=1/2$, a compromise between stability and accuracy is achieved, where the condition number still increases but not as drastically as in the case $\beta=0$, and the relative error decreases as the mesh is refined.

\begin{table}[h]
\centering
\captionsetup{font=small}
\caption {\label{tab:table2} Condition number \eqref{Eqn.CondNumber} and relative error \eqref{Eqn.RelError} for mesh refinements $\EPW = 12$, $24$, $48$ and $96$ using the IGA discretization, and for various choices for the parameter $\beta$ to govern the virtual sources displacement distance $h \sim \lambda / \EPW^{\beta}$. The wavenumber $k=4 \pi$ is fixed and $\lambda = 2 \pi / k$ is the wavelength. The number of Gauss nodes per element $\NG=4$ is also fixed.} 
\begin{tabular}{lllllllll}
\hline
 & \multicolumn{2}{c}{$\beta=0$} & & \multicolumn{2}{c}{$\beta=1/2$} & &  \multicolumn{2}{c}{$\beta=1$}  \\ 
\cline{2-3} \cline{5-6} \cline{8-9} 
$\EPW$ &  Cond. Number  & Rel. Error  &  & Cond. Number  & Rel. Error &  &  Cond. Number  & Rel. Error  \\  

$12$   &  $6.14 \times 10^{0}$  & $1.76 \times 10^{-5}$  & & $6.14 \times 10^{0}$  & $1.76 \times 10^{-5}$ & & $6.14 \times 10^{0}$  & $1.76 \times 10^{-5}$  \\

$24$   &  $2.74 \times 10^{2}$  & $2.33 \times 10^{-7}$  & & $5.19 \times 10^{1}$ & $6.78 \times 10^{-7}$   & & $1.68 \times 10^{1}$  & $1.25 \times 10^{-5}$  \\
  
$48$   &  $6.13 \times 10^{5}$  & $3.53 \times 10^{-13}$ & & $8.76 \times 10^{2}$ & $6.72 \times 10^{-9}$  & & $3.96 \times 10^{1}$  & $1.14 \times 10^{-5}$  \\

$96$   &  $2.86 \times 10^{12}$  & $6.57 \times 10^{-15}$ & & $3.31 \times 10^{4}$ & $5.18 \times 10^{-11}$  & & $6.68 \times 10^{1}$  & $1.56 \times 10^{-5}$  \\
  \hline
\end{tabular}
\end{table}

\begin{figure}[h]
\centering
\captionsetup{font=small}
\subfloat[$k=4 \pi$]{
\includegraphics[height=0.24 \textwidth, trim = 2 0 2 0, clip]{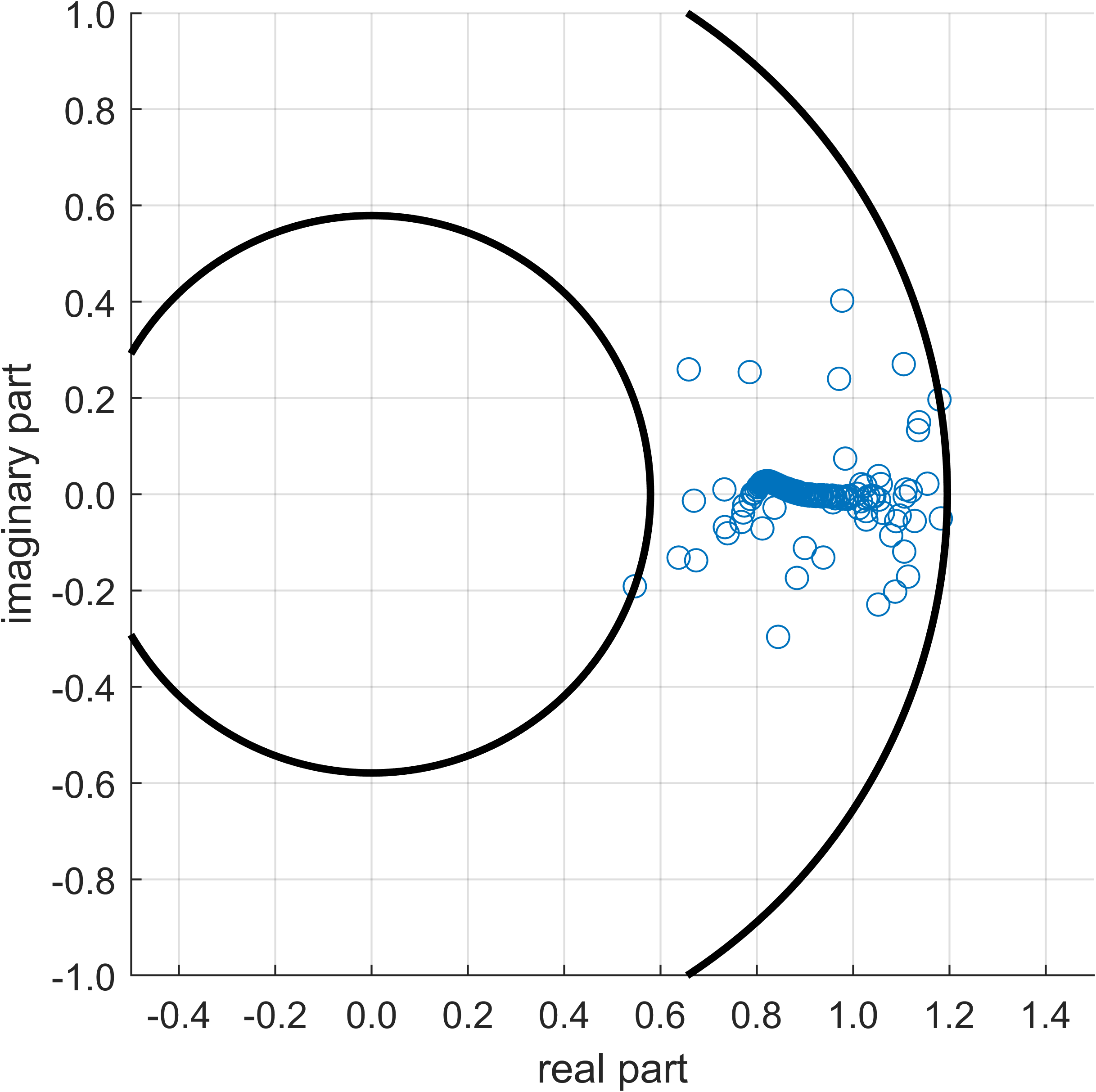}
\label{subfig:K4pi2}
}
\subfloat[$k=8 \pi$]{
\includegraphics[height=0.24 \textwidth, trim = 2 0 2 0, clip]{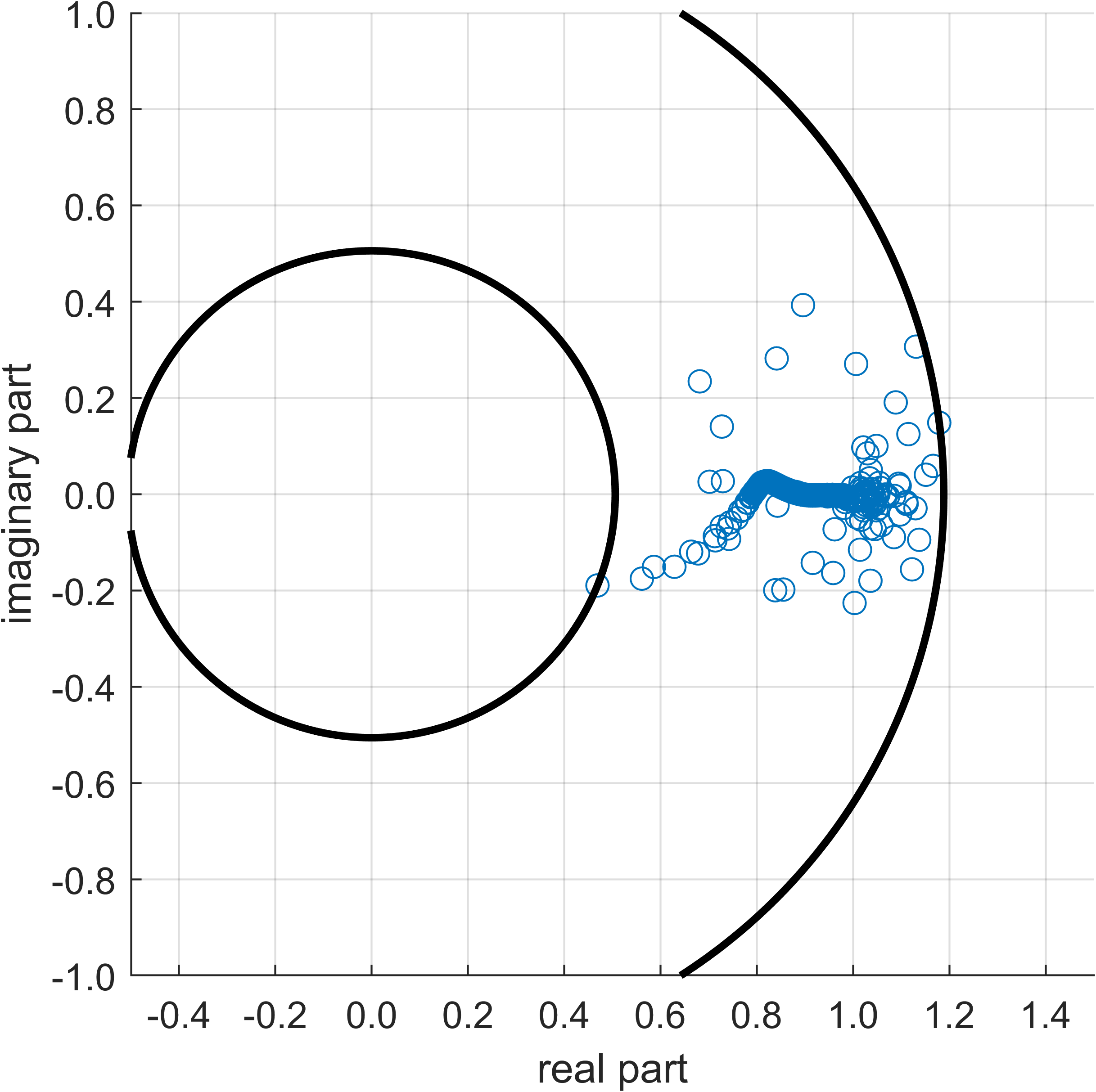}
\label{subfig:K8pi2}
}
\subfloat[$k=16 \pi$]{
\includegraphics[height=0.24 \textwidth, trim = 2 0 2 0, clip]{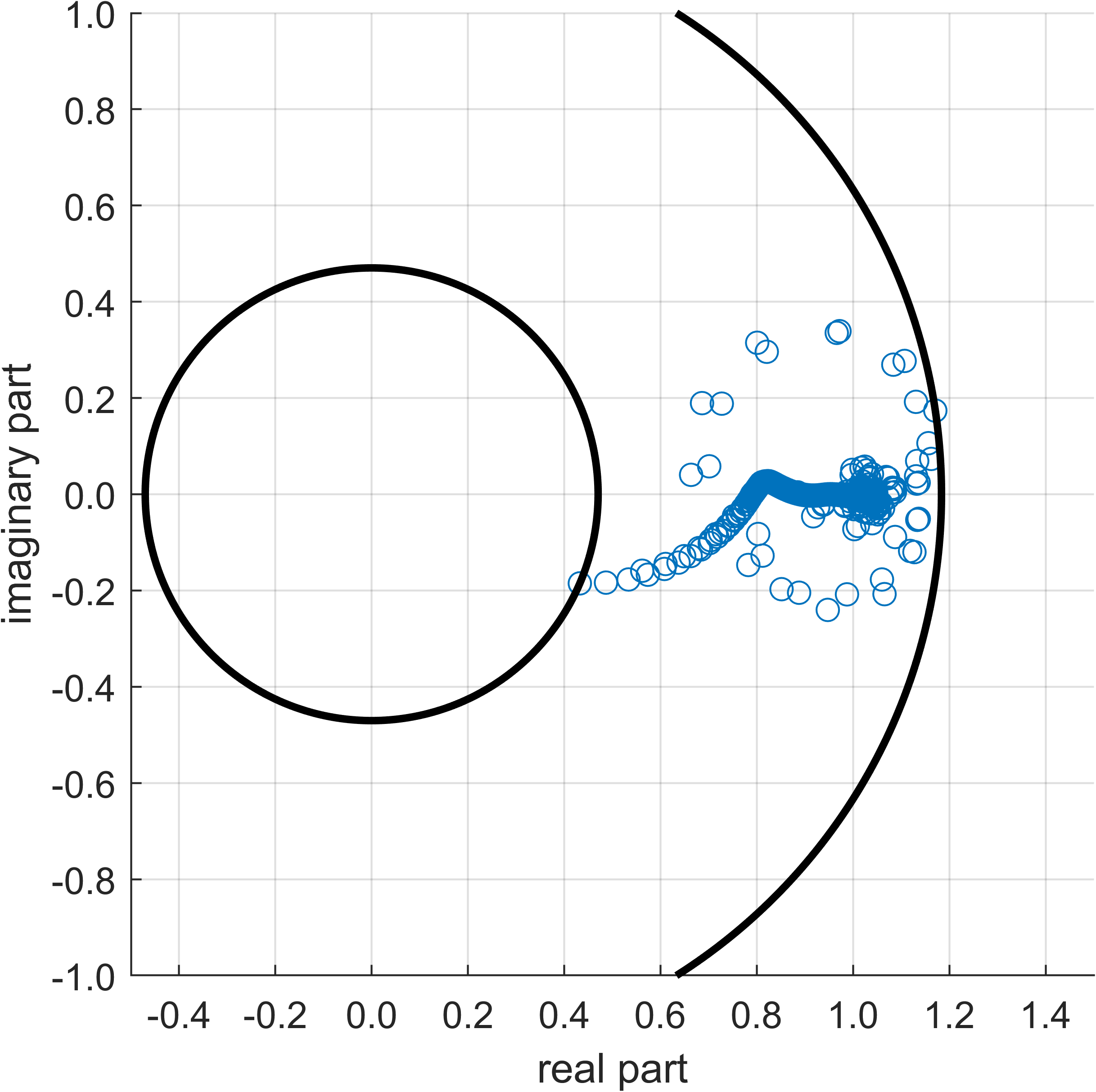}
\label{subfig:K16pi2}
} 
\subfloat[$k=32 \pi$]{
\includegraphics[height=0.24 \textwidth, trim = 2 0 2 0, clip]{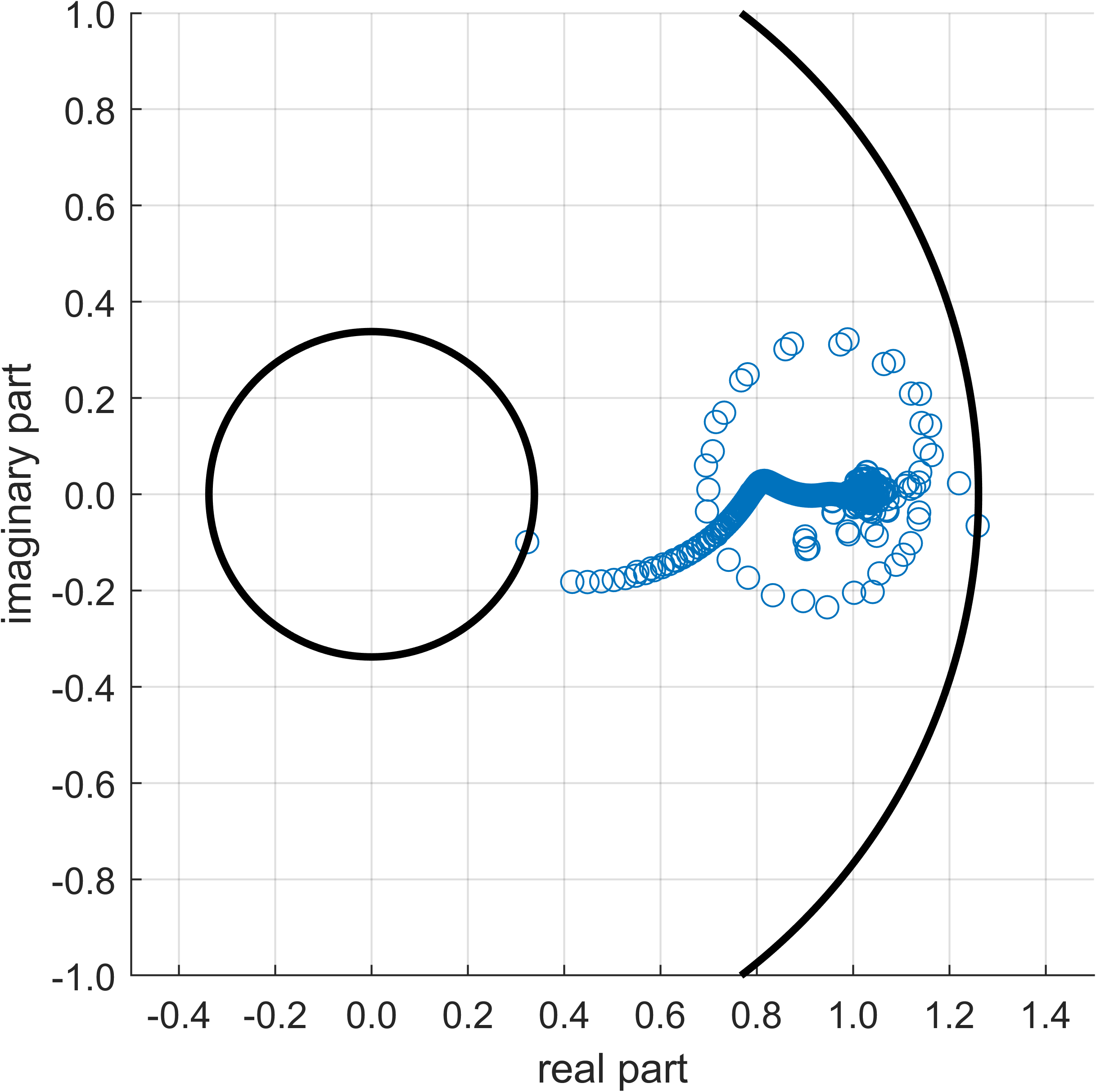}
\label{subfig:K32pi2}
}
\caption{Spectrum of the IGA implementation of the matrix $\BB_{N}\AA_{N}$ for wavenumbers $k=4 \pi$, $k=8 \pi$, $k=16 \pi$, and $k=32 \pi$. The circles' radii correspond to the largest and smallest magnitude of the spectrum. 
In all cases, the discretization was done using $\EPW = 12$ along the boundary $\Gamma$. The corresponding degrees of freedom are $N=150$, $N=302$, $N=604$, and $N=1206$, respectively. The virtual source displacement distance is given by $h = \lambda / 12$ where $\lambda = 2\pi / k$ is the wavelength.}
\label{fig:IGASpectrumIncreasingFreq}
\end{figure}

The spectrum of the proposed matrix $(\BB_{N} \AA_{N})$ obtained after the IGA discretization is shown in Figure \ref{fig:IGASpectrumIncreasingFreq} for increasing wavenumber $k=4\pi$, $8\pi$, $16\pi$, and $32\pi$ but fixed number of elements per wavelength $\EPW=12$. The displacement distance for the virtual sources is $h=\lambda/12$ where $\lambda=2\pi/k$ is the wavelength. As in the previous section, here we also observe that the upper- and lower-bounds of the spectrum are practically independent of the frequency. As analyzed in the Appendix \ref{Section.GaussQuadPhi}, we expect the numerical quadrature to remain accurate near the shifted singularity while keeping the condition number of the matrix $(\BB_{N} \AA_{N})$ relatively small for this choice of $h$ as the wavenumber $k$ increases, as long as $\EPW$ is kept constant. This is an advantage over the conventional method of virtual sources which has 
eigenvalues approach zero exponentially fast as $k$ grows. See \cite{Kress1986} or the discussion in the Appendix \ref{Section.Spectral}.

\subsection{Scattering of a plane wave from a circular obstacle}
\label{Sec.SubsectionPlane_IGA}

As a second numerical experiment, we compute the numerical solution for a canonical wave scattering problem with an exact analytical solution. We consider a plane wave 
\begin{align}
    \uinc(x,y) = e^{i k x} \label{Eqn.PlaneWave}
\end{align}
propagating along the $x$-axis. Using the notation from Section \ref{Sec.MathFormulation}, this wave impinges on an acoustically soft impenetrable obstacle $\Omega^{-} \subset \mathbb{R}^2$ with circular boundary $\Gamma$ of radius $R=1$. The total field $\utot = \uinc + \usc = 0$ on $\Gamma$ and satisfies the Helmholtz equation in $\Omega^{+}$ so that the scattered field $\usc$ solves \eqref{Eqn.001}-\eqref{Eqn.003} with $f = - \uinc$ on $\Gamma$. The exact solution for this problem can be found for instance in \cite[Ch. 4]{Martin-Book-2006}. For this numerical experiment, we are interested in the behavior of the error as both the number of elements per wavelength $\EPW$ and the number of Gauss quadrature nodes per element $\NG$ increase. The relative error for these numerical runs are displayed in Table \ref{tab:tablePlaneWave}. We observe that for a fixed $\EPW$, the error decreases initially as $\NG$ increases and then it stagnates. This phenomenon is consistent with the ability of the Gauss quadrature to integrate accurately the profile of the Green's function with the displaced singularity. Once this quadrature error is diminished, then the error introduced by the finite-dimensional NURBS basis dominates. This dominant error can decreases by refining the mesh, that is, increasing $\EPW$. Overall, this proposed method of virtual sources provides excellent accuracy for this canonical benchmark.

\begin{figure}[h]
\centering
\captionsetup{font=small}
\subfloat[Numerical solution]{
\includegraphics[height=0.28 \textwidth, trim = 0 0 0 0, clip]{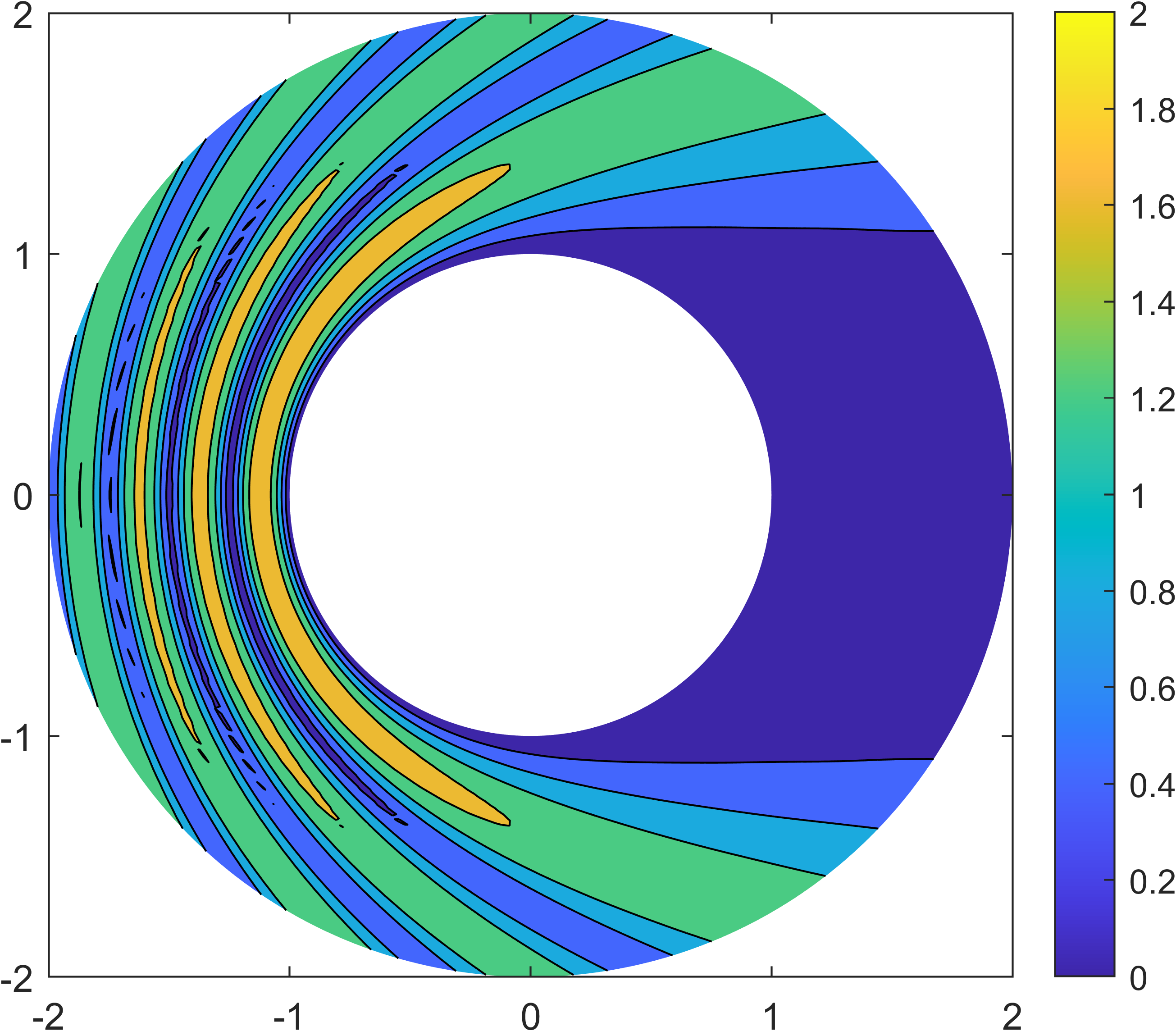}
\label{subfig:NumSol3}
}
\subfloat[Error ($\log_{10}$ scale) ]{
\includegraphics[height=0.28 \textwidth, trim = 0 0 0 0, clip]{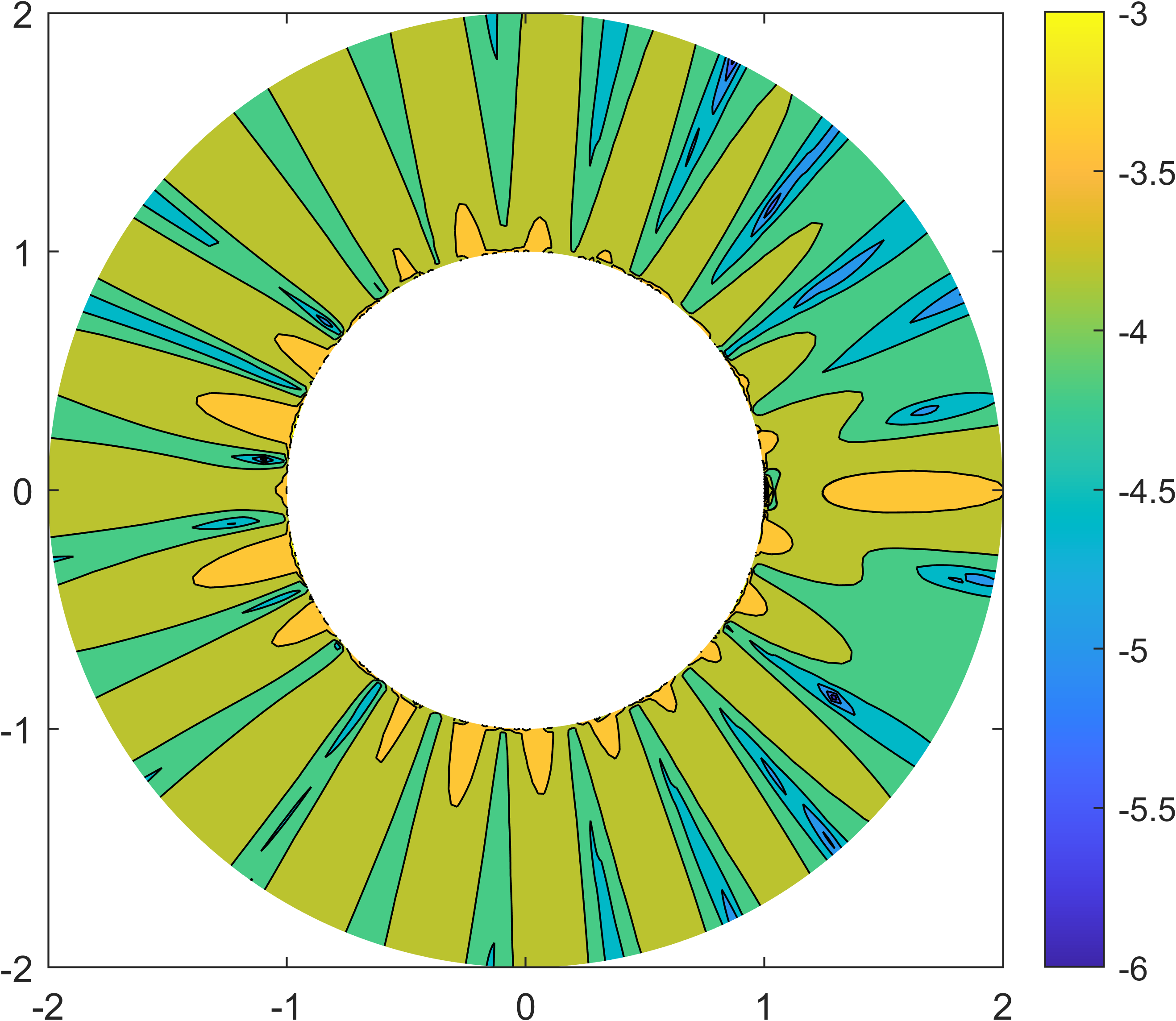}
\label{subfig:Error3}
}
\subfloat[GMRES residual]{
\includegraphics[height=0.28 \textwidth, trim = 0 0 0 0, clip]{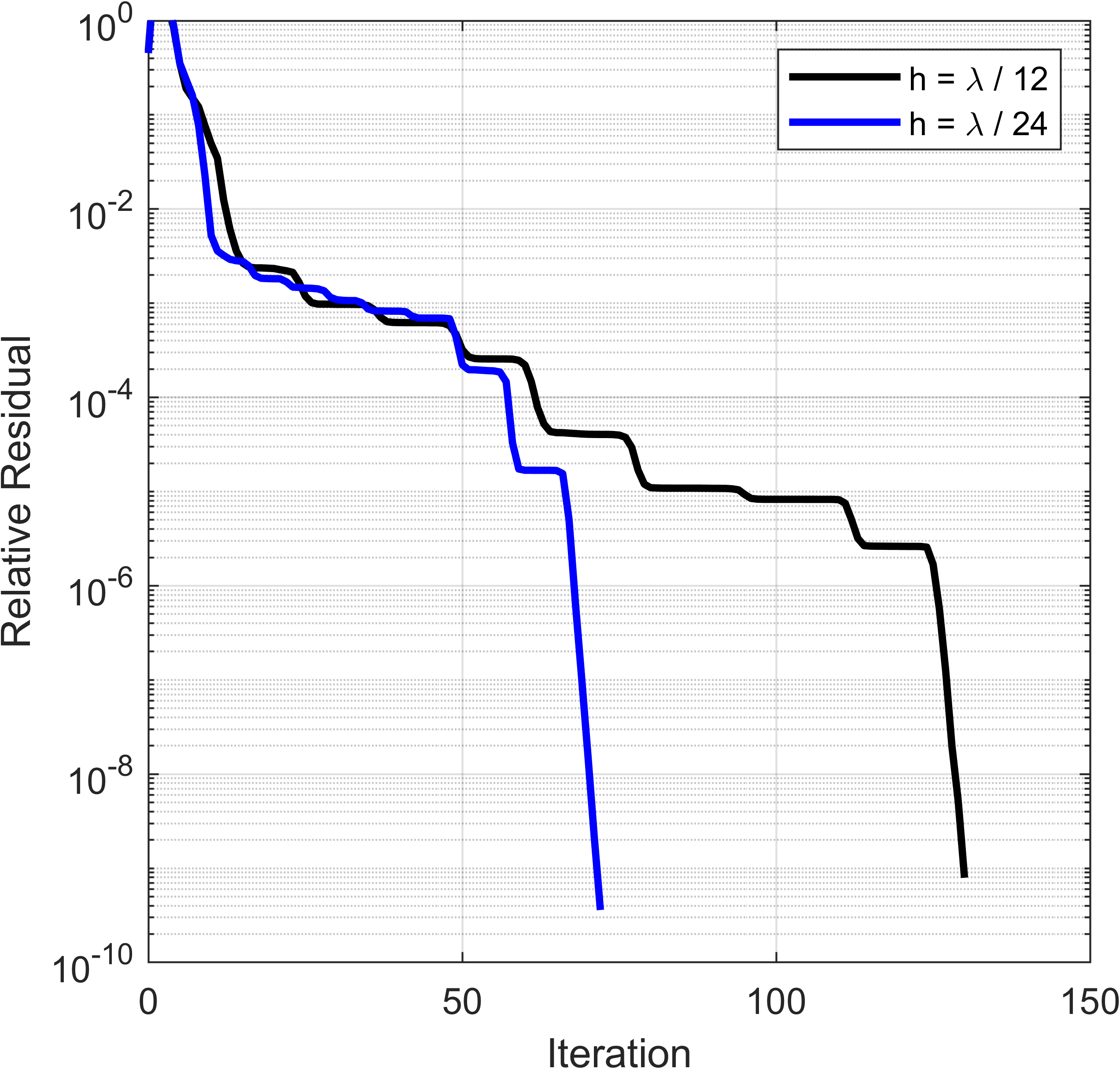}
\label{subfig:GMRES_res3}
}
\caption{ (a) Numerical solution (total field) using the IGA implementation for $k=4\pi$, (b) the error ($\log_{10}$ scale) in a neighborhood of the boundary $\Gamma$, and (c) the residual for the GMRES iterations for the cases $h = \lambda/12$ and $h = \lambda/24$. The numerical solution was obtained using the proposed method of virtual sources with $\EPW=12$, a Pad\'e approximation $\Losrc$ of the DtN map using $4$ terms, and Gaussian quadrature with $\NG = 4$ nodes per elements.}
\label{Fig.IGA_ExactSol_and_comparison_Plane}
\end{figure}

\begin{table}[h]
\centering
\captionsetup{font=small}
\caption { \label{tab:tablePlaneWave} Relative error for the scattering of a plane wave from a circular obstacle described in subsection \ref{Sec.SubsectionPlane_IGA} for various values of the number of elements per wavelength $\EPW$ and number of Gauss nodes per element $\NG$.} 
\begin{tabular}{cccc}
\hline
$\NG$ & $\EPW=6$ & $\EPW=12$ &  $\EPW=24$ \\  

$4$ & $ 6.23\e{-3} $ & $ 2.44\e{-4} $ & $ 7.60\e{-8} $  \\

$5$ & $ 2.23\e{-3} $ & $ 1.19\e{-5} $ & $ 3.49\e{-8} $  \\
  
$6$ & $ 6.86\e{-4} $ & $ 1.93\e{-6} $ & $ 1.80\e{-8} $  \\

$7$ & $ 2.36\e{-4} $ & $ 2.97\e{-7} $ & $ 1.58\e{-8} $  \\

$8$ & $ 8.37\e{-5} $ & $ 1.76\e{-7} $ & $ 1.36\e{-8} $  \\

\hline
\end{tabular}
\end{table}

\subsection{Scattering of a plane wave from a 3D obstacle}

\begin{figure}[h]
\centering
\captionsetup{font=small}
\includegraphics[height=0.45 \textwidth, trim = 0 0 0 0, clip]{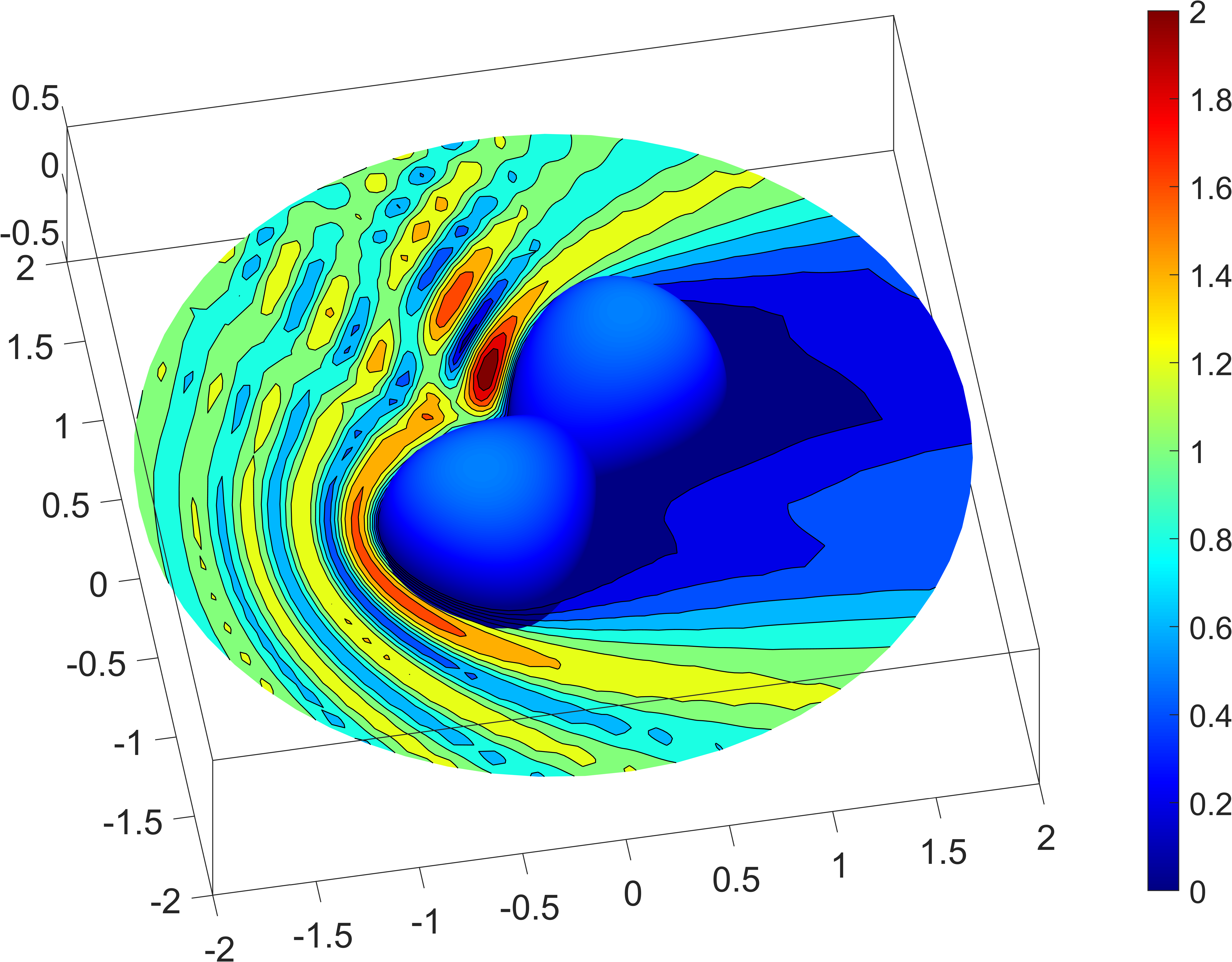}
\caption{Total field on the $xy$-plane for three-dimensional wave scattering from an acoustically soft obstacle.}
\label{fig:3D_scattering}
\end{figure}

Finally, we also present the result of a numerical run for a three-dimensional scenario. We consider an obstacle whose boundary is a surface of revolution. This is shown in Figure \ref{fig:3D_scattering}. The incident wave field is a plane wave with wavenumber $k=2 \pi$. The numerical solution was computed using the IGA discretization with NURBS basis of degree 2. The obstacle's surface was meshed using $\EPW = 6$ elements per wavelength approximately, which led to $N=1845$
degrees of freedom for the scattering problem. As before, the OSRC approximation $\Losrc$ to the DtN map $\Lambda$ was computed using Pad\'e approximations with $4$ terms as detailed in the Appendix \ref{Section.Pade}. The virtual source displacement distance $h = \lambda / 24$ where $\lambda = 2 \pi / k$ is the wavelength.  Upon discretization, equation \eqref{Eqn.106} was numerically solved using GMRES which took less than 100 iterations for the residual to drop below $10^{-10}$. 
This run demonstrates the feasibility of the proposed virtual source method to compute solution to wave scattering problems in three-dimensional settings.

\section{Conclusion}
\label{Section.Conclusion}

In this paper we proposed a novel method of virtual sources to formulate boundary integral equations with continuous kernels for exterior wave propagation problems. Single- and double-layer potentials were combined through the use of OSRC to obtained well-conditioned systems upon discretization. However, by contrast to classical boundary integral formulations, 
we displaced the singularity of the Green's function by a small distance $h>0$ to a virtual parallel surface. As a result, the discretization was be performed on the actual physical boundary with continuous kernels.

The virtual displacement parameter $h>0$ controls the conditioning of the discrete system. We provided mathematical understanding to choose $h$, in terms of the wavelength and mesh refinements, in order to achieve a balance between accuracy and stability. We implemented the proposed method of virtual sources using piecewise linear elements and IGA formulations, observed exceptionally well-behaved spectra, and solved the corresponding systems using matrix-free GMRES iterations. The results were compared to analytical solutions for canonical problems. We conclude that the proposed method leads to accurate solutions and robust conditioning/stability properties for a wide range of wavelengths and mesh refinements.

\appendix

\section{Appendices}

\subsection{Pad\'e approximation of the DtN map}
\label{Section.Pade}
Here we offer some details about Pad\'e approximation of the pseudodifferential symbol of the DtN map. Recall that the principal symbol of the DtN map is given by $ik \sqrt{1 - \sigma^2/k^2}$ where $-\sigma^2$ represent the pseudodifferential symbol of the surface Laplacian. Hence, we seek the Pad\'e approximant for the function $(1+z)^{1/2}$ about $z=0$. In partial fraction form, the Pad\'e approximant is
\begin{align} \label{Eqn.PadeAppendix}
P_{M}(z) =  a_{0} + \sum_{m=1}^{M} \frac{a_{m} }{ z - b_{m} }
\end{align}
where the real-valued coefficients $a_{m}$ and $b_{m}$ are computed in order to match the value of the function $(1+z)^{1/2}$ and its first $2M$ derivatives at $z=0$ \cite{BakerBook1996}. Unfortunately, these Pad\'e approximants suffer from  instabilities when $z \leq -1$ due to the branch cut along the negative real line starting at $z=-1$ \cite{Lu1998a,Lu1998b,Milinazzo1997,Modave2020}. Milinazzo et al. \cite{Milinazzo1997} avoided this problem by considering the Pad\'e approximation of $e^{i \theta/2 } \left( 1 + \zeta \right)^{1/2}$ where $\zeta = (1+z) e^{-i \theta} - 1$, which has a rotated branch cut defined by the angle $\theta$. As a result, the approximation remains stable and continuous for $z \in \mathbb{R}$.  This $\theta$-rotated Pad\'e approximant of order $M$, in partial fraction form, is given by
\begin{align} \label{Eqn.PadeRotatedAppendix}
P_{M,\theta}(z) = e^{i \theta/2} \left( a_{0} + \sum_{m=1}^{M} \frac{a_{m} }{ (1+z) e^{- i \theta} - 1 - b_{m} } \right) = \tilde{a}_{0} + \sum_{m=1}^{M} \frac{\tilde{a}_{m} }{ z - \tilde{b}_{m} }
\end{align}
where the complex-valued Pad\'e coefficients are defined by $\tilde{a}_{0} = a_{0} e^{i \theta /2 }$, $\tilde{a}_{m} = a_{m} e^{3i \theta /2}$ and $\tilde{b}_{m} = (1+b_{m})e^{i \theta} - 1$ for $m=1,2,...,M$ \cite{BakerBook1996,Milinazzo1997}.


\subsection{Numerical quadrature near the shifted singularity of the fundamental solution}
\label{Section.GaussQuadPhi}

Here we explore some of the details concerning the choice of the source displacement distance $h >0$ for the proposed method of virtual sources. We assume that the error $E$ associated with approximating the integral of a generic smooth function $g$ over an interval $[a,b]$ is bounded as follows
\begin{align}
E \leq C_{n} (b-a)^{2n+1} \max_{y \in [a,b]} | g^{(2n)}(y) | \label{Eqn.GenericGaussQuad}
\end{align}
for some constant $C_{n}$, where the superscript $^{(2n)}$ denotes the $2n$-th derivative, and where $n$ is the number of quadrature nodes. This is the case for Gaussian quadrature or the trapezoidal rule for smooth periodic functions \cite[\S 12.1]{Kress-Book-1999}. 
For our purposes, the integrand takes the form
\begin{align}
g(y) = \partial_{\nu(y)} \Phi(x,y - h\nu(y)) v(y) - \Phi(x,y - h\nu(y)) \Losrc v(y)
\end{align}
where $x$ is taken as a fixed parameter. Assuming that the function $v$ oscillates like a wave field with wavenumber $k$, then the $m$-th derivative of $v$ behaves like $v^{(m)} \sim k^m v$. Since the operator $\Losrc$ is of degree 1, then we can estimate that $(\Losrc v)^{(m)} \sim k^{m+1} v$. For the shifted kernel of the single-layer potential, we have the following behavior near the singularity $\Phi \sim \ln(kh)$ and $\Phi^{(m)} \sim h^{-m}$ for $m\geq 1$. For the shifted kernel of the double-layer potential, we have $(\partial_{\nu(y)} \Phi) ^{(m)} \sim h^{-(m+1)}$. As a result, we have that
\begin{align}
|g^{(2n)}(y)| \lesssim \sum_{m=0}^{2n} D_{m} \left( \frac{k^{2n-m}}{h^{m+1}} + \frac{k^{2n -m + 1}}{h^m} \right)
\end{align}
for some constants $D_{m}$. Here the symbol $\lesssim$ is being used to denote inequality up to a multiplicative constant independent of $g$, $k$ or $h$. 
Now we assume that the size of an element is given in terms of the chosen number of elements per wavelength ($\EPW$), 
\begin{align}
(b-a) = \frac{\lambda}{\EPW}, \qquad \text{where $\lambda = 2 \pi / k$ is the wavelength.}
\end{align}
Hence, the quadrature error near the shifted singularity is bounded as follows
\begin{align}
E \lesssim \frac{C_{n}}{\EPW^{2n+1}} \sum_{m=0}^{2n} D_{m} \left(  \frac{\lambda}{h} + 1 \right) \left( \frac{\lambda}{h} \right)^m.
\end{align}
Therefore, in order to maintain the order of accuracy for the  quadrature with $n$ nodes near the shifted singularity, it is sufficient that the shifting distance for the displacement of the virtual sources be in the order of the wavelength of the oscillatory fields,
\begin{align}
h \sim \lambda. \label{Eqn.lambda_h}
\end{align}


\subsection{Spectral characterization of the method of virtual sources}
\label{Section.Spectral}

Recall that the governing operator behaves as in \eqref{Eqn.102}. Also let the DtN map have an eigenvalue decomposition with eigenvalues $\{ \gamma_{n} \}_{n=1}^{\infty}$. Since $\Lambda$ is a pseudo-differential operator of degree $1$ whose principal symbol is $i \sqrt{k^2 - \xi^2}$, then we expect its eigenvalues $\{ \gamma_{n} \}$ to grow linearly as $n \to \infty$. For the propagating modes, the eigenvalues $|\gamma_{n}| < k$ and are purely imaginary. For the evanescent modes, the eigenvalues $| \gamma_{n} | > k$ and possess a negative real part. As a result, the governing operator $\AA \approx \exp ( h \Lambda )$ has eigenvalues that approach zero exponentially fast, that is $\exp ( h \text{Re}(\gamma_{n}) )$, as $n \to \infty$. This is the ill-conditioning introduced by displacing the virtual sources a distance $h > 0$ away from the physical boundary $\Gamma$.

For numerical implementations, we consider a discrete approximation  $\AA_{N}$ of $\AA$ employing a corresponding discrete approximation $\Lambda_{N}$ of $\Lambda$ that correctly captures the behavior of the first $N$ eigenvalues of $\Lambda$. Here $N$ represents the degrees of freedom introduced in Sections \ref{Section.LinearDiscrete} and \ref{Section.IGADiscrete}, which is proportional to the number of elements per wavelength $\EPW$.
This means that for $N$ sufficiently large, $\AA_{N}$ has an eigenvalue close to zero,  $\sim \exp (- C h N )$ for some constant $C>0$. Hence, we could control the ill-conditioning of the discrete version of the method of virtual sources by choosing $h$ small enough at the expense of sacrificing some quadrature accuracy. This leads to replacing \eqref{Eqn.lambda_h} by
\begin{align}
h \sim \lambda/N^{\beta} \sim \lambda / \EPW^{\beta}
\label{Eqn.lambda_h_N}
\end{align}
for some $0 \leq \beta \leq 1$. In the extreme case $\beta=0$, we expect the numerical quadrature to be accurate but the governing system to be ill-conditioned. At the other extreme $\beta=1$, we expect the numerical quadrature to lose most of its accuracy in favor of  well-conditioning of the scheme.

These arguments can be made precise for a circular or spherical boundary. We assume that $\Gamma$ is a circle of radius $a$, and the virtual parallel surface $\Gamma_{h}$ if circle of radius $a-h$. A radiating solution $u$ to the Helmholtz equation in the exterior of $\Gamma_{h}$ can be expressed in Fourier modes as follows,
\begin{align}
u(r,\theta) =  \sum_{n=-\infty}^{\infty} v_{n} \frac{H_{|n|}(kr)}{H_{|n|}(k(a-h))} e^{i n \theta}, \qquad \text{where} \quad v_{n} = \frac{1}{2 \pi} \int_{0}^{2\pi} v(\theta) e^{-i n \vartheta} d \vartheta,
\end{align}
where $H_{n}$ denotes the Hankel function of the first kind of order $n$, and $v$ is the density of sources on $\Gamma_{h}$. However, as explained in Section \ref{Sec.MVCS}, since $\Gamma_{h}$ and $\Gamma$ are sufficiently close parallel surfaces, then $v$ can be identified as a source density on $\Gamma$. Hence, the governing operator can be seen as $\AA : H^{0}(\Gamma) \to H^{0}(\Gamma)$ given by
\begin{align}
(\AA v)(\theta) =  \sum_{n=-\infty}^{\infty} v_{n} \frac{H_{|n|}(ka)}{H_{|n|}(k(a-h))} e^{i n \theta}
\end{align}
when using the exact DtN map $\Lambda$ in the definition \eqref{Eqn.IntOp}. The above reveals the spectral behavior of the operator $\AA$. Its eigenvalues are then given by
\begin{align}
\lambda_{n} = \frac{H_{|n|}(ka)}{H_{|n|}(k(a-h))}, \qquad n=0, \pm 1, \pm 2 , ... .
\end{align}
We are interested in the spectral behavior for large eigenvalues. For $|n| \gg |z|^2$, the Hankel functions behave like
\begin{align}
H_{|n|}(z) \sim \frac{1}{|n| !} \left( \frac{z}{2} \right)^{|n|} - i \frac{(|n|-1)!}{\pi} \left( \frac{2}{z} \right)^{|n|}
\end{align}
which implies that the eigenvalues behave asymptotically as 
\begin{align}
\lambda_{n} \sim \left( 1 - \frac{h}{a} \right)^{|n|}, \qquad \text{as $|n| \to \infty$}.
\end{align}
This expression confirms the exponential decay of the eigenvalues which governs the ill-conditioning of the method of virtual sources for $h>0$. For a numerical method that effectively truncates the eigenvalue decomposition to $N$ terms, then choosing $h$ in the form \eqref{Eqn.lambda_h_N}
would prevent the exponential decay of the first $N$ eigenvalues. We note, however, that slower decay (such as polynomial) can still be present.


\vfill
\clearpage

\bibliographystyle{plain}
\bibliography{library}


\end{document}